\input vanilla.sty
\scaletype{\magstep1}
\scalelinespacing{\magstep1}
\def\bull{\vrule height .9ex width .8ex depth -.1ex}


\title Complex interpolation of Hardy-type subspaces\\
\endtitle

\author N.J.  Kalton\footnote{Research supported by NSF
grant DMS-9201357\newline AMS Classifications:  46B70,
46E30} \\ Department of Mathematics\\ University of
Missouri-Columbia\\ Columbia, Mo. 65211 \endauthor

\subheading{Abstract} We consider the problem of complex
interpolation of certain Hardy-type subspaces of K\"othe
function spaces.  For example, suppose $X_0$ and $X_1$ are
K\"othe function spaces on the unit circle $\bold T,$ and
let $H_{X_0}$ and $H_{X_1}$ be the corresponding Hardy
spaces.  Under mild conditions on $X_0,X_1$ we give a
necessary and sufficient condition for the complex
interpolation space $[H_{X_0},H_{X_1}]_{\theta}$ to coincide
with $H_{X_{\theta}}$ where $X_{\theta}=[X_0,X_1]_{\theta}.$
We develop a very general framework for such results and our
methods apply to many more general sitauations including the
vector-valued case.

\vskip2truecm

\subheading{1.  Introduction} Let $X$ be a K\"othe function
space on the circle $\bold T$ equipped with its usual Haar
measure.  Consider the Hardy subspace $H_X$ consisting of
all $f\in X\cap N^+$ where $N^+$ is the Smirnov class or
Hardy algebra.  Provided $X\subset L_{\log}$ (see Section 2
for the definition) this is a closed subspace.  Consider the
following two problems:\newline (1) When is $H_X$
complemented in $X$ by the usual Riesz projection?\newline
(2) If $X_0,X_1$ are two such K\"othe function spaces when
is it true that the complex interpolation space
$X_{\theta}=[X_0,X_1]_{\theta}$ satisfies
$H_{X_{\theta}}=[H_{X_0},H_{X_1}]_{\theta}?$

In the case of weighted $L_p-$spaces, a precise answer to
(1) was given by Muckenhaupt [26] in terms of the so-called
$A_p-$conditions.  In the case $p=2,$ the Helson-Szeg\"o
theorem [15] gives an alternative precise criterion; in the
same direction Cotlar and Sadosky [8], [9] gave necessary
and sufficient conditions for all weighted $L_p-$spaces (see
also [10]).  Subsequently, Rubio de Francia extended the
Cotlar-Sadosky methods to all 2-convex or 2-concave K\"othe
function spaces.  In the case of $L_p$-spaces (without
weights) (2) is answered by a well-known theorem of Jones
[16], [17] (cf. recent proofs by Xu [34], M\"uller [27] and
Pisier [29]).  For weighted $L_p-$spaces (2) has recently
been studied by Cwikel, McCarthy and Wolff [11] and
Kisliakov and Xu [20], [21] (who also consider vector-valued
analogues).  See also [33].

In this paper we will develop a very general approach to
question (2) by relating it to (1).  We will be able to
answer (2) completely under some mild restrictions on the
spaces.  In fact our approach uses very little specific
information about Hardy spaces or properties of analytic
functions and we give our results in a rather general
setting, which includes abstract Hardy spaces generated by
weak$^*$-Dirichlet algebras and certain vector-valued cases
as studied by Kisliakov and Xu.

We limit our discussion in the introduction to the case of
the circle.  Let us say that a K\"othe function space $X$ is
{\it BMO-regular} if and only if there exist constants
$(C,M)$ so that given $0\le f\in X$ there exists $g\ge f$
with $\|g\|_X \le M\|f\|_X$ and $\|\log g\|_{BMO}\le C.$ A
weighted $L_p-$space, $L_p(w)$ is BMO-regular if and only if
$\log w\in BMO.$ The concept of BMO-regularity appears
implicitly first in the work of Cotlar and Sadosky [9] and
also in Rubio de Francia [30] in connection with the
boundedness of the Hilbert transform (it should be noted
that in both cases the boundedness of the Hilbert transform
is related to the BMO-regular of a space derived from $X$,
not of $X$ itself).  We show that a super-reflexive space
$X$ is BMO-regular if and only if the Riesz projection is
bounded on an interpolation space $L_2^{\theta}X^{1-\theta}$
for some $\theta>0$ (cf.  [18] for other conditions
equivalent to this property for $X$).

If $X_0,X_1$ are super-reflexive and
$X_0,X_1,X_0^*,X_1^*\subset L_{\log}$ then we give a
necessary and sufficient condition for
$H_{X_{\theta}}=[H_{X_0},H_{X_1}]_{\theta}$ where
$0<\theta<1$ (in this case we say that the Hardy algebra
$H=N^+$ is interpolation stable at $\theta$ for
$(X_0,X_1)$).  Consider first the case when $X_0$ is
BMO-regular; then it is necessary and sufficient that $X_1$
is BMO-regular.  For the general case the necessary and
sufficient condition is obtained by ``lifting'' the
direction $X_0\to X_1$ to create a parallel direction
$L_2\to Z;$ the condition is then that $Z$ is BMO-regular.
This is precisely explained in Section 5; let us remark that
if $X_1=wX_0$ is obtained by a change of weight then
$Z=w^{1/2}L_2=L_2(w^{-1})$ so that the condition is simply
that $\log w\in BMO$.  Our result includes the results of
the previous work of Kisliakov and Xu [21] as special cases
and extends, as we have explained, to a very general
setting, thus giving also vector-valued applications.

Let us also comment on the methods used.  In Section 3 we
discuss a very general formulation of question (2):  when
does the operation of interpolation commute with taking a
particular subspace?  Our main result is that if this
happens, then under appropriate conditions, one can
extrapolate the boundedness of a projection onto the
subspace.  In Section 4 we consider an arbitrary
self-adjoint operator $T$ on $L_2.$ We then discuss for
which K\"othe spaces $X$ it is true that $T$ is bounded on
$L_2^{1-\theta}X^{\theta}$ for some $\theta>0.$ If we assume
that $T$ is bounded on some $L_p$ where $p\neq 2$ then this
can be answered in terms of the weighted $L_2-$spaces on
which $T$ is bounded.  These results are of course closely
related to the earlier work of Cotlar and Sadosky [9] and
Rubio de Francia [30]; unlike [30] we do not assume
2-convexity or 2-concavity but our conclusions are somewhat
weaker.

We put these ideas together in Section 5, restricting our
attention to ``Hardy-type'' algebras, which we introduce as
an abstraction of the Smirnov class; in this case our
operator $T$ becomes the orthogonal projection onto $H_2.$
We are then able to relate the results of Section 4 to the
notion of BMO-regularity and prove our main results.  We
discuss further applications in Section 6. At the end of
Section 6, we improve the results of Kisliakov and Xu ([20],
[21]) on interpolation of vector-valued Hardy spaces, by
giving necessary and sufficient conditions for such
interpolation to be ``stable'' at least in the
super-reflexive case.

Let us mention that we use some ideas from [18]; however we
have tried to avoid using differential techniques in order
to keep our approach as simple as possible.  We plan a
further paper showing how by using such techniques one can
improve and extend these results.  We do use freely however,
the notion of an indicator function for a K\"othe function
space as introduced and studied in [18].

We would like to thank Michael Cwikel, Mario Milman and
Richard Rochberg for discussing this problem with us, and
Cora Sadosky for some very helpful remarks.  We also thank
Quanhua Xu for several important suggestions which we have
incorporated.

\vskip2truecm

\subheading{2.  K\"othe function spaces} Let $S$ be a Polish
space and let $\mu$ be a probability measure on $S.$ Let
$L_0(\mu)$ denote, as usual, the space of all equivalence
classes of (complex) Borel functions on $S$ with the
topology on convergence in measure.

We define a K\"othe quasinorm on $L_0$ to be a
lower-semicontinuous functional $f\to \|f\|_X$ defined on
$L_0$ with values in $[0,\infty]$ such that:\newline (1)
$\|f\|_X=0$ if and only if $f=0$ a.e.\newline (2)
$\|f\|_X\le \|g\|_X$ whenever $|f|\le |g|$ a.e.\newline (3)
There exists a constant $C$ so that $\|f+g\|_X\le
C(\|f\|_X+\|g\|_X)$ for $f,g\in L_0.$\newline (4) There
exists $u\in L_0$ so that $u>0$ a.e. and
$\|u\|_X<\infty.$\newline Associated to the K\"othe
quasi-norm we can associate a maximal quasi-K\"othe function
space $X=\{f:\|f\|_X<\infty\}.$ $X$ is then a quasi-Banach
space under quasi-norm $f\to \|f\|_X;$ furthermore
$B_X=\{f:\|f\|_X\le 1\}$ is closed in $L_0$ so that $X$ has
the Fatou property (cf.  [24]).  We can also define a
minimal quasi-K\"othe function space $X_0$ to be the closure
of $L_{\infty}\cap X$ in $X.$ In this paper, however, we
will only deal with maximal spaces (i.e. spaces with the
Fatou property).  If, in (3), $C=1$ then $B_X$ is convex and
$X$ is a Banach space; in this case we say that $X$ is a
(maximal) K\"othe function space.  {\it Henceforward we will
adopt the convention that all spaces are maximal.}

If $X$ is a K\"othe function space and $w\in L_{0,\bold R}$
with $w>0$ a.e. we define the weighted space $wX$ by
$\|f\|_{wX}=\|fw^{-1}\|_X.$ Thus $wL_p=L_p(w^{-p}).$

If $X$ is a K\"othe function space we will let $X^*$ denote
its K\"othe dual i.e. the maximal K\"othe function space
induced by $\|\,\|_{X^*}$ where $\|f\|_{X^*}=\sup_{g\in
B_X}\int |fg|d\mu.$ It is not difficult to show that $X^*$
is also a K\"othe function space.  Of course, if $X$ is
reflexive as will usually be the case, then $X^*$ is the
Banach dual of $X.$

We recall that a quasi-K\"othe function space $X$ is
$p$-convex where $0<p<\infty$ with constant $M$ if for every
$f_1,\ldots,f_n\in X$ we have:  $$ \|(\sum_{k=1}^n
|f_k|^p)^{1/p}\|_X \le M(\sum_{k=1}^n\|f_k\|_X^p)^{1/p},$$
and $q$-concave $(0<q<\infty)$ with constant $M$ if for
every $f_1,\ldots,f_n\in X,$ $$
(\sum_{k=1}^n\|f_k\|_X^q)^{1/q} \le M \|(\sum_{k=1}^n
|f_k|^q)^{1/q}\|_X.$$ If $X$ is $p$-convex and $q$-concave
there is an equivalent quasi-norm so that the $p$-convexity
and $q$-concavity constants are both one.  For conevenience
we will say that $X$ is {\it exactly} $p$-convex or
$q$-concave if the associated constant of convexity or
concavity is one.  $X$ is a K\"othe function space if and
only it is 1-convex with constant one.  A K\"othe function
space is super-reflexive if and only if it is $p$-convex and
$q$-concave for some $1<p\le q<\infty.$

For any K\"othe function space $X$ we define the
quasi-K\"othe function space $X^{\alpha}$ by
$\|f\|_{X^{\alpha}}= \||f|^{1/\alpha}\|_X^{\alpha}.$ Then
$X^{\alpha}$ is exactly $1/\alpha-$convex.  If $X,Y$ are two
K\"othe function spaces and $0< \alpha,\beta <\infty$ we can
define a quasi-K\"othe function space
$Z=X^{\alpha}Y^{\beta}$ by setting
$\|f\|_Z=\inf\{\max(\|g\|_X,\|h\|_Y)^{\alpha+\beta}:
|f|=|g|^{\alpha}|h|^{\beta}\}.$ Then $Z$ is exactly
$1/(\alpha+\beta)-$convex.  It may also be shown easily
that, since $X,Y$ are assumed maximal, there is always an
optimal factorization $|f|=|g|^{\alpha}|h|^{\beta}.$

We now describe a simple method of doing calculations with
K\"othe function spaces introduced in [18].  We will not
need the full force of the results in [18] and thus we will
try to keep the description brief.  Let us recall [18] that
a semi-ideal $\Cal I$ is a subset of $L_{1,+}$ so that if
$0\le f\le g\in \Cal I$ then $f\in \Cal I$; $\Cal I$ is
strict if it contains a strictly positive function.  For a
functional $\Phi:\Cal I\to \bold R$ we define
$\Delta_{\Phi}(f,g)=\Phi(f)+\Phi(g)-\Phi(f+g).$ We say that
$\Phi$ is semilinear if:  \newline (1) If $f\in \Cal I,$ and
$\alpha>0$ then $\Phi(\alpha f)=\alpha \Phi(f),$\newline (2)
There is a constant $\delta$ so that for all $f,g\in\Cal I$
we have $\Delta_{\Phi}(f,g)\le
\delta(\|f\|_1+\|g\|_1).$\newline (3) If $f\in \Cal I$ and
$0\le f_n\le f$ with $\|f_n\|_1\to 0$ then
$\lim\Phi(f_n)=0.$

If $X$ is a K\"othe function space we define $\Cal I_X$ to
be the set of nonnegative functions $f$ in $L_1$ so that
$\sup_{x\in B_X}\int f\log_+|x|d\mu<\infty$ and there exists
$x\in B_X$ so that $f\log |x|$ is integrable.  Then $\Cal
I_X$ is a strict semi-ideal.

On $\Cal I_X$ we can define the indicator functional
$\Phi_X=\sup_{x\in B_X}\int f\log|x|\,d\mu.$ The indicator
function $\Phi_X$ is semilinear with $\delta\le \log 2$ (see
[18], Proposition 4.2).  In the special case $X=L_1$ we
obtain $\Cal I_X=L\log L$ and $\Phi_{L_1}(f)=\Lambda(f)=\int
f\log fd\mu.$ It then may be shown that for general $X$ and
$f,g\in \Cal I_X\cap L\log L,$ we have $0 \le
\Delta_{\Phi_X}(f,g)\le \Delta_{\Lambda}(f,g).$

There is a converse to this result ([18], Theorem 5.2).  If
$\Phi$ is a semilinear map on a strict semi-ideal $\Cal
I\subset L\log L$ so that $0\le \Delta_{\Phi}(f,g)\le
\Delta_{\Lambda}(f,g)$ for all $f,g\in \Cal I$ then there is
a unique K\"othe function space $X$ so that $\Cal I_Z\supset
\Cal I$ and $\Phi(f)=\Phi_X(f)$ for $f\in \Cal I.$
Furthermore $X$ is exactly $p$-convex and exactly
$q$-concave if and only if $\frac1q\Delta_{\Lambda}(f,g)\le
\Delta_{\Phi_X}(f,g)\le \frac1p\Delta_{\Lambda}(f,g)$ for
$f,g\in\Cal I.$

It is also easy to see that if $Z=X^{\alpha}Y^{\beta}$ then
$\Phi_Z(f)=\alpha\Phi_X(f)+\beta\Phi_Y(f)$ for $f\in \Cal
I_X\cap\Cal I_Y.$ This enables us to use the indicator
functions to calculate spaces.  Let us give a simple
application which we will use later; basically this is a
simple generalization of Pisier's extrapolation theorem [28]
(cf.  [18], Corollary 5.4).

\proclaim{Proposition 2.1}Suppose $X,Y$ are K\"othe function
spaces with $Y$ super-reflexive.  Then there exists a
super-reflexive K\"othe function space $Z$ and $0<\theta<1$
so that $Y=X^{1-\theta}Z^{\theta}$ up to equivalence of
norm.\endproclaim

\demo{Proof}For convenience we do all calculations on a
strict semi-ideal contained in $\Cal I_X\cap \Cal I_Y\cap
L\log L.$ We can assume that $Y$ is exactly $p$-convex and
$q$-concave where $\frac1p+\frac1q=1$ and $2<q<\infty.$ Let
$\alpha=\frac{1}{2q}.$ We define
$\Phi=\Phi_Y+\alpha(\Phi_Y-\Phi_X).$ Then $$\Delta_{\Phi}=
(1+\alpha)\Delta_{\Phi_Y}-\alpha\Delta_{\Phi_X} \ge
\alpha\Delta_{\Lambda}.$$ Similarly, $$\Delta_{\Phi}\le
(1+\alpha)(1-2\alpha)\Delta_{\Lambda}\le
(1-\alpha)\Delta_{\Lambda}.$$ Thus we can apply Theorem 5.2
of [18] to find a space $Z$ so that $\Phi_Z=\Phi$ and $Z$
will be $2q-$concave and $r$-convex where
$\frac1r+\frac1{2q}=1.$\enddemo

Finally let us define $L_{\log}$ to be the Orlicz space of
all $f\in L_0$ so that $\int \log_+|f|d\mu<\infty.$ Then
$L_{\log}$ can be F-normed by $f\to \int \log(1+|f|)d\mu.$
We will be especially concerned with the class of K\"othe
function spaces $\Cal X$ of all $X$ so that $X,X^*\subset
L_{\log}.$

\proclaim{Lemma 2.2}If $X$ is a K\"othe function space, then
the following conditions are equivalent:\newline (1)
$X\in\Cal X.$\newline (2) If $f\in X$ there exists $g\in X$
with $g\ge |f|,$ and $\log g\in L_1.$\newline (3) If
$\epsilon>0$ and $f\in X$ there exists $g\in X$ with $g\ge
|f|$, $\|g\|_X\le \|f\|_X+\epsilon$ and $\log g\in L_1.$ (4)
$L_{\infty}\subset \Cal I_X.$\endproclaim

\demo{Proof}$(1)\Rightarrow (4):$ We must show $\chi_S\in
\Cal I_X.$ If $X\subset L_{\log}$ then it follows from the
Closed Graph Theorem that the inclusion is continuous and
hence $\sup_{x\in B_X}\int f\log_+|x|d\mu<\infty.$ On the
other hand by a theorem of Lozanovskii [25] (cf.  [13],
[25]) there exist nonnegative $x\in B_X$ and $x^*\in
B_{X^*}$ with $xx^*=\chi_S.$ Thus $\log
|x|=\log_+|x|-\log_+|x^*|\in L_1.$

$(4)\Rightarrow (3):$ If $f\in X$ then $\log_+|f|\in L_1.$
However there exists $h\in B_X$ with $\log |h|\in L_1.$ Now
take $g=\max(|f|,\eta |h|)$ for small enough $\eta.$

$(3)\Rightarrow (2):$ Obvious.

$(2)\rightarrow (1):$ Clearly the conditions imply $X\subset
L_{\log}.$ Now suppose $x^*\in X^*.$ There exists $x\in X$
with $\log |x|\in L_1.$ Now $xx^*\in L_1$ so that $\int\log
|xx^*|d\mu<\infty.$ Hence $\int \log |x^*|d\mu<\infty$ i.e.
$\log_+|x^*|\in L_1.$ Thus $X^*\subset
L_{\log}.$\bull\enddemo

\demo{Remark}In doing calculations with indicator functions
we can always restrict to a small enough strict semi-ideal.
Later in the paper for economy we will not mention the
domains of the indicators in question when doing algebraic
manipulations.  The reader may wish simply to consider only
spaces $X\in\Cal X$ and regard all indicator functions as
defined on $L_{\infty,\bold R}.$\enddemo

\vskip2truecm

\subheading{3.  Complex interpolation of subspaces} Let us
describe a very general setting for complex-type
interpolation.  We recall that if $X$ is a topological
vector space and $D$ is an open subset of the complex plane
then a function $F:D\to X$ is analytic if for each $a\in D$
there exists a neighborhood $U$ of $a$ and a power series
$\sum_{n=0}^{\infty} x_n(z-a)^n$ so that $F(z)=\sum
x_n(z-a)^n$ for $z\in U.$ We will consider a triple
$(D,X,\Cal F)$ where $D$ is an open subset of the complex
plane conformally equivalent to the unit disk, $X$ is a
complex topological vector space and $\Cal F$ is subspace of
the space $\Cal A(D,X)$ of all $X$-valued analytic functions
on $D$ equipped with a norm $F\to \|F\|_{\Cal F}$ such that:
\item{1.} If $F\in \Cal A(D,X)$ and $\varphi$ is any
conformal mapping of $D$ onto $\Delta$ then $\varphi F\in
\Cal F$ if and only if $F\in \Cal F$ and $\|\varphi
F\|_{\Cal F}=\|F\|_{\Cal F}.$ \item{2.} If $z\in D$ and
$x\in X$ then $\inf\{\|F\|_{\Cal F}:\ F(z)=x\}=0$ if and
only if $x=0.$

Under these assumptions, we define, for $z\in D$,
$X_z=\{x:\|x\|_{X_z}<\infty\}$ where
$\|x\|_{X_z}=\inf\{\|F\|_{\Cal F}:\ F(z)=x\}.$ We will call
the spaces $\{X_z:z\in D\}$ the {\it interpolation field}
generated by $\{D,X,\Cal F\}.$

The following elementary lemma will be used repeatedly.  If
$a,b\in D$ we let $\delta(a,b)=|\varphi(b)|$ where $\phi$ is
any conformal map of $D$ onto $\Delta$ with $\varphi(a)=0.$
Thus if $D=\Delta$ we have
$$\delta(a,b)=\frac{|b-a|}{|1-\bar ab|}.$$ If $D=\Cal S$ is
the strip $\Cal S=\{z:0<\Re z<1\}$ then for $0\le s,t\le 1$
we have $$\delta(s,t) =\frac{\sin \frac{\pi}2|s-t|}{\sin
\frac{\pi}2(s+t)}.$$

\proclaim{Lemma 3.1}Suppose $F,G\in \Cal F$ and that $a\in
D.$ Suppose $F(a)=G(a).$ Then $\|F(z)-G(z)\|_{X_z} \le
\delta(a,z)\|F-G\|_{\Cal F}$.\endproclaim

\demo{Proof}Let $\varphi$ be a conformal map of $D$ onto
$\Delta$ with $\varphi(a)=0.$ Let $H\in \Cal F$ be defined
so that $\varphi H=F-G.$ Then $\|H\|_{\Cal F}\le
\|F-G\|_{\Cal F}$ and $\|F(z)-G(z)\|_{X_z} \le
|\varphi(z)|\|H\|_{\Cal F}.$ The lemma follows.\bull\enddemo

Now suppose that $V$ is a linear subspace of $X.$ Let $\Cal
F(V)$ be the space of $F\in \Cal F$ such that $F(z)\in V$
for every $z\in D.$ Let $(V_z)$ be the interpolation field
generated by $\{D,Y,\Cal F(V)\}.$ We will say that $V$ is
{\it interpolation stable} at $z\in D$ if there is a
constant $C$ so that for $v\in V_z$ we have $\|v\|_{V_z}\le
C\|v\|_{X_z}.$ The least such constant $C$ we denote by
$K(z)=K(z,V)$ where $K(z)=\infty$ if $V$ fails to be
interpolation stable.

\proclaim{Theorem 3.2}Suppose $V$ is interpolation stable at
some $a\in D;$ let $K(a)=K.$ Then $V$ is interpolation
stable at any $z\in D,$ with $3K\delta(a,z)< 1,$ and $K(z)
\le 4K(1-3K\delta)^{-1}.$ In particular, the set of $z\in D$
so that $V$ is interpolation stable at $z$ is
open.\endproclaim

\demo{Proof} Suppose $\delta=\delta(a,z)<1/(3K).$ Suppose
$v\in V_z$; then for $\epsilon>0$ we can pick $F\in \Cal
F(V)$ and $G\in \Cal F$ so that $F(z)=G(z)=v$ and
$\|F\|_{\Cal F}\le (1+\epsilon)\|v\|_{V_z}$ while
$\|G\|_{\Cal F}\le (1+\epsilon)\|v\|_{X_z}.$ Thus by Lemma
3.1, $$ \|F(a)-G(a)\|_{X_a} \le
(1+\epsilon)\delta(\|v\|_{X_z}+\|v\|_{V_z}).$$ It follows
that $$ \|F(a)\|_{X_a} \le
(1+\epsilon)((1+\delta)\|v\|_{X_z} +\delta \|v\|_{V_z}).$$

Now pick $H\in \Cal F(V)$ so that $H(a)=F(a)$ and
$\|H\|_{\Cal F}\le (1+\epsilon)\|F(a)\|_{V_a} \le
(1+\epsilon)K \|F(a)\|_{X_a}.$ Then $$ \align \|v\|_{V_z}
&\le \|F(z)-H(z)\|_{V_z} + \|H(z)\|_{V_z}\\ &\le
\delta(\|F\|_{\Cal F} +\|H\|_{\Cal F}) + \|H\|_{\Cal F}\\
&\le \delta \|F\|_{\Cal F} + (1+\delta)(1+\epsilon)K
\|F(a)\|_{X_a}\\ &\le
\delta(1+K+K\delta)(1+\epsilon)\|v\|_{V_z}
+(1+\delta)^2(1+\epsilon)^2K\|v\|_{X_z}.  \endalign $$

Since $\delta<1$ and $K\ge 1,$ we have $$ \|v\|_{V_z} \le
3K\delta\|v\|_{V_z} + 4K \|v\|_{X_z}$$ whence we conclude
that $$ \|v\|_{V_z} \le \frac{4K}{(1-3K\delta)}\|v\|_{X_z}$$
and so $K(z) \le 4K(1-3K\delta)^{-1}.$ \bull\enddemo

\proclaim{Theorem 3.3}Suppose $V,W$ are two subspaces of $X$
which are both interpolation stable at $a.$ Suppose further
that $X_a=V_a\oplus W_a.$ Then there exists $\eta>0$ so that
if $\delta(a,z)<\eta$ then $X_z=V_z\oplus W_z.$\endproclaim

\demo{Proof} Let $K=\max(K(a,V),K(a,W))$ and let
$M=\max(\|P\|,\|Q\|)$ where $P,Q$ are the induced
projections from $X_a$ onto $V_a$ and $W_a.$ We let
$\eta=1/(300K^2M).$ Suppose $z\in D$ satisfies
$\delta=\delta(a,z)<\eta.$

First suppose $z$ is such that $V_z+W_z$ fails to be dense
in $X_z.$ Then there exists $x\in X_z$ so that
$\|x\|_{X_z}=1$ and $\|x-(v+w)\|_{X_z}\ge \frac12$ whenever
$v\in V_z$ and $w\in W_z.$ Pick any $F\in\Cal F$ with
$F(z)=x$ and $\|F\|_{\Cal F}\le 2.$ Then there exist
$G\in\Cal F(V)$ and $H\in\Cal F(W)$ so that $\|G\|_{\Cal
F},\|H\|_{\Cal F}\le 4KM$ and $G(a)=PF(a),\ H(a)=QF(a).$
Then $$ \align \frac12 &\le \|F(z)-G(z)-H(z)\|_{X_z}\\ &\le
\delta(2 + 2(4KM))\\ &\le 10KM\eta.  \endalign $$ This
contradiction immediately leads to the conclusion that
$V_z+W_z$ is dense in $X_z.$

To complete the proof, suppose $v\in V_z,\ w\in W_z$
satisfying $\|v+w\|_{X_z}=1.$ Let
$\gamma=\max(1,\|v\|_{X_z},\|w\|_{X_z}).$ We will show that
$\gamma\le 8KM$ and this will complete the proof.  We choose
$F\in \Cal F$ with $\|F\|_{\Cal F}\le 2$ and $F(z)=v+w.$
Notice that $K(z,V),K(z,W)\le 8K$ by Theorem 3.3.  We
therefore pick $G\in \Cal F(V),\ H\in \Cal F(W)$ so that
$\|G\|_{\Cal F}\le 10K\gamma$ and $\|H\|_{\Cal F}\le
10K\gamma,$ and $G(z)=v,\ H(z)=w.$

Now we have, by Lemma 3.1, $\|F(a)-G(a)-H(a)\|_{X_a} \le
30K\gamma\delta$ (where $\delta=\delta(a,z)$) and hence
$\|PF(a)-G(a)\|_{X_a} \le 30 KM\gamma\delta$ and
$\|QF(a)-H(a)\|_{X_a} \le 30 KM\gamma\delta.$ However
$\|PF(a)\|_{X_a}\le M\|F(a)\|_{X_a} \le 2M$ and so we obtain
an estimate $$ \|G(a)\|_{X_a} \le 2M+ 30KM\gamma\delta)$$
with a similar estimate for $\|H(a)\|_{X_a}.$ Thus $$
\|G(a)\|_{V_a} \le 2KM+30K^2M\gamma\delta< 2KM
+\frac{\gamma}{10}$$ and we can find $E\in\Cal F(V)$ with
$E(a)=G(a)$ and $\|E\|_{\Cal F}\le 3KM +
\frac{1}{10}\gamma.$ Now $\|v\|_{X_z} \le \|E(z)\|_{X_z}
+\|E(z)-G(z)\|_{X_z}.$ Thus $$ \align \|v\|_{X_z} &\le
(1+\delta)\|E\|_{\Cal F} +\delta \|G\|_{\Cal F}\\ &\le
3KM(1+\delta) + \frac{1}{5}\gamma + 10K\gamma\delta \\ &\le
4KM + \frac12 \gamma.  \endalign $$ With a similar estimate
on $w$ we obtain $$ \gamma \le (4KM) + \frac12\gamma$$ and
so $\gamma\le 8KM$ as promised.\bull\enddemo

Let us now give a simple application.  Obviously one special
case of the above construction is the usual Calder\`on
method of complex interpolation.  To be more precise, note
that if $(X_0,X_1)$ is a Banach couple then if we take
$D=\Cal S,$ $X=X_0+X_1$ and $\Cal F$ to be the space of
functions $F\in\Cal A(S,X)$ so that $F$ is bounded on $\Cal
S$ and extends continuously to the closure of $\Cal S$ so
that $F$ is $X_j$-continuous on the line $\Re z=j$ for
$j=0,1$ then the interpolation field generated given by
$X_z=[X_0,X_1]_{\theta}$ where $\theta=\Re z.$

\proclaim{Theorem 3.4}Suppose $(X_0,X_1)$ is a Banach
couple.  Suppose for some $0<\theta<1,$ we have that
$[X_0,X_1]_{\theta}$ is a UMD-space and
$[H_2(X_0),H_2(X_1)]_{\theta}=H_2[X_0,X_1]_{\theta}.$ Then
there exists $\eta>0$ so that if $|\phi-\theta|<\eta$ then
$[X_0,X_1]_{\phi}$ is also UMD.\endproclaim

\demo{Remark}Blasco and Xu [2] show if $X_0$ and $X_1$ are
UMD-spaces then
$[H_2(X_0),H_2(X_1)]_{\theta}=H_2(X_{\theta}).$ This result
is therefore a converse to their result.  They also present
an example of Pisier to show that
$[H_2(X_0),H_2(X_1)]_{\theta}$ need not coincide with
$H_2(X_{\theta}).$ We remark that in [18] we construct an
example where $[X_0,X_1]_{1/2}=L_2$ but $X_{\theta}$ is not
UMD for any $\theta\neq\frac12,$ thus giving another
counterexample.\enddemo

\demo{Proof}We consider the Banach couple
$(L_2(X_0),L_2(X_1)).$ Let $V$ be the subspace of
$L_2(X_0)+L_2(X_1)\subset L_2(X_0+X_1)$ of all $f$ so that
$$ \int_0^{2\pi} e^{int}f(e^{it}) dt =0$$ for $n>0$.  Let
$W$ be the space of all $f$ so that $$ \int_0^{2\pi}
e^{int}f(e^{it}) dt =0$$ for $n\le 0.$ Our assumptions
guarantee that $V,W$ are interpolation stable at $\theta$
and that
$$[L_2(X_0),L_2(X_1)]_{\theta}=L_2([X_0,X_1]_{\theta}=V_{\theta}\oplus
W_{\theta}.$$ By Theorem 3.3 we obtain a similar
decomposition for $|\phi-\theta|<\eta$ which implies the
result.\bull\enddemo

Let us now discuss the case of K\"othe function spaces.
Suppose $S$ is a Polish space and $\mu$ is a probability
measure on $S.$ As in [18] we consider the class $\Cal N^+$
of all functions $F:\Delta\to L_0$ of the form
$F(z)(s)=F_s(z)$ where $F_s$ is in the Smirnov class $N^+$
for almost every $s\in S.$ Then $\Cal N^+(\Cal S)$ is the
class of maps $F:\Cal S\to L_0$ where $F\circ\varphi\in \Cal
N^+$ where $\varphi:  D\to\Cal S$ is any conformal mapping.
If $F\in\Cal N^+( \Cal S)$ we can extend $F$ to the lines
$z=j+it$ ($j=0,1$) so that $F(j+it)=\lim_{s\to j}F(s+it)$ in
$L_0,$ for a.e.  $t.$

Suppose $X_0,X_1$ are K\"othe function spaces (assumed
maximal so that $f\to \|f\|_{X_j}$ is lower-semi-continuous
on $L_0).$ Consider the space $\Cal F=\Cal F(X_0,X_1)$ of
all $F\in\Cal N^+(\Cal S)$ so that $\|F\|_{\Cal
F}=\max_{j=0,1}\{\text{ess }\sup\|F(j+it)\|_{X_j}\}<\infty.$
Then $\Cal F$ generates an interpolation field $X_z$ for
$z\in\Cal S$ so that $X_z=X_0^{1-\theta}X^{\theta}$ where
$\theta=\Re z.$ Now suppose $Z$ is a separable K\"othe
function space which contains both $X_0$ and $X_1$.  It is
essentially shown in [18] that if $F\in\Cal F(X_0,X_1)$ then
$F:\Cal S\to Z$ is analytic and $\lim_{s\to
j}F(s+it)=F(j+it)$ in the space $Z$ (so that we can work in
$Z$ in place of $L_0$); see Lemma 2.2 of [18].

Now suppose $V$ is a linear subspace of $L_0$ so that for
some separable K\"othe function space $Z\supset X_0,X_1$ the
space $V\cap Z$ is closed in $Z.$ Then $V_j=V\cap X_j$ is
closed in $X_j$ for $j=0,1.$ Furthermore if $F\in \Cal
F(X_0,X_1;V)=\Cal F(X_0,X_1)(V)$ then $F$ has boundary
values in $V_j$ along the line $z=j+it,\ -\infty<t<\infty.$
The method of interpolation generated this way is not
precisely the complex method introduced by Calder\'on, but
we now make some remarks which establish that under
reasonable hypotheses we obtain the same result.

The usual interpolation spaces $[V_0,V_1]_{\theta}$ are
induced by considering the subspace $\Cal F_c(X_0,X_1;V)$ of
all $F$ so that (a) $F$ is analytic into $V_0+V_1,$ (b)
$\lim_{s\to j}F(s+it)$ exist a.e. in $V_0+V_1$ and (c) $t\to
F(j+it)$ is Bochner measurable into $V_j$ for $j=0,1.$ In
fact only condition (c) is required; this is a consequence
of the following lemma.

\proclaim{Lemma 3.5}If $G\in \Cal F(X_0,X_1;V)$ then $G\in
\Cal F_c(X_0,X_1;V)$ if and only if for each $j$ the map
$t\to G(j+it)$ has essentially separable range in $X_j.$
\endproclaim

\demo{Proof}This is essentially proved in Lemma 2.2 of [18],
although our assumptions are a little less strict; we sketch
the argument.  Observe first that $t\to G(j+it)$ is Bochner
measurable into $X_j$ for each $j.$ Let $\varphi:\Delta
\to\Cal S$ be a conformal mapping.  Suppose $G\in\Cal
F(X_0,X_1;L_0).$ Then $F=G\circ\varphi\in\Cal N^+$ and has
$L_0-$boundary values $F(e^{it})$ for a.e.  $0\le t<2\pi.$
Then $t\to F(e^{it})$ is $V_0+V_1$ and $X_0+X_1$-measurable.
It is thus Bochner integrable in both $X=X_0+X_1,$ and
$W=V_0+V_1.$

Now suppose $w\in X^*$ is strictly positive.  It follows
that $$ \int_0^{2\pi}\int_S
|F_s(e^{it})|w(s)d\mu(s)\frac{dt}{2\pi} <\infty.$$ Then if
$F_s(z)=\sum_{n\ge 0}a_n(s)z^n$ for $|z|<1$ we observe that
$t\to F_s(e^{it})\in N^+\cap L_1=H_1$ for a.e.  $s.$ Hence
if $n\ge 0,$ we have, $\mu-$a.e., $$ \int_0^{2\pi}
F_s(e^{it})e^{-int}\frac{dt}{2\pi} = a_n(s).$$ Similarly if
$n<0$ we have $$ \int_0^{2\pi}
F_s(e^{it})e^{-int}\frac{dt}{2\pi} = 0.$$ Hence we can also
evaluate the Bochner integrals in $L_1(w)$ $$ \int_0^{2\pi}
F(e^{it})e^{-int}\frac{dt}{2\pi} = a_n$$ when $n\ge 0$ and
$$ \int_0^{2\pi} F(e^{it})e^{-int}\frac{dt}{2\pi} = 0$$ when
$n<0.$ These integrals have the same values in $W$ and $X$
and it thus follows easily that $F:\Delta\to W$ is analytic
and has a.e. boundary values $F(e^{it}).$ This implies that
$G\in \Cal F_c(X_0,X_1;V).$\bull\enddemo

\demo{Remark}It follows that if $X_0,X_1$ are separable
K\"othe function spaces (which are as usual assumed to have
the Fatou property) then
$[X_0,X_1]_{\theta}=X_0^{1-\theta}X_1^{\theta}.$ (cf.
[5.]).

\proclaim{Proposition 3.6} Suppose $X_0,X_1$ are K\"othe
function spaces and that $V$ is a linear subspace of $L_0$
so that $V\cap Z$ is closed in $Z$ for some separable
K\"othe function space $Z$ containing $X_0,X_1.$ Suppose
$0<\theta<1.$ Suppose either \newline (a) $X_0$ and $X_1$
are both separable or \newline (b)
$X_{\theta}=X_0^{1-\theta}X_1^{\theta}$ is
reflexive.\newline Then $V$ is interpolation stable at
$\theta$ for the interpolation method generated by $\Cal
F(X_0,X_1)$ if and only if $[V_0,V_1]_{\theta}=V\cap
X_{\theta}$ up to equivalence of norm.  \endproclaim

\demo{Proof}It is immediately clear that
$[V_0,V_1]_{\theta}=V\cap X_{\theta}$ implies the
interpolation stability of $V.$ In the other direction
consider first (a).  In this case Lemma 3.5 implies that
$\Cal F_c(X_0,X_1;V)=\Cal F(X_0,X_1;V)$ and the conclusion
is immediate.

Now consider (b); let $V_{\theta}$ be the space induced by
the method $\Cal F(X_0,X_1;V)$ and let
$W=[V_0,V_1]_{\theta}.$ Let $B$ be the closed unit ball of
$V_{\theta}$ and let $B'$ be the closed unit ball of $W.$
Then $B'\subset B;$ it follows from the Open Mapping Theorem
since both spaces are complete that if we can show that $B'$
is dense in $B$ then $W$ and $V_{\theta}$ coincide.  Suppose
$B''$ is the $V_{\theta}-$closure of $B'.$ Since
$V_{\theta}$ is a closed subspace of $X_{\theta},$ $B''$ is
weakly compact.  If $B''\neq B$ there exists $v\in
V_{\theta}$ with $\|v\|_{V_{\theta}}\le 1$ and $v\notin
B''.$ Since $B''$ is weakly compact in $Z$ there exists
$\phi\in Z^*$ so that $|\phi(v)|>1$ but $\sup_{w\in
B'}|\phi(w)|\le 1.$ Let $|\phi(v)|=\rho^2$ where $\rho>1.$
Pick any $F\in\Cal F(X_0,X_1;V)$ with $\|F\|_{\Cal F}\le
\rho$ and $F(\theta)=v.$

Now for $\tau>0$ let $$G_{\tau}(z)=\frac1\tau
\int_0^{\tau}F(z+it)\,dt,$$ for $z\in\Cal S,$ where the
integrals are computed in $Z.$ It is clear that the boundary
values of $G_{\tau}$ are given by the same formula; hence
$t\to G_{\tau}(j+it)$ is continuous in $X_j$ for $j=0,1.$
Since clearly $G_{\tau}\in \Cal F(X_0,X_1;V)$, Lemma 3.5 can
be applied to give that $G_{\tau}\in\Cal F_c(X_0,X_1;V).$
Hence $\rho^{-1}G_{\tau}(\theta)\in B'.$ We conclude that
$|\phi(G_{\tau}(\theta))| \le \rho$ and letting $\tau\to 0$
gives $|\phi(v)|\le \rho$, a contradiction.\bull\enddemo

In the situations when we will apply this result we will
consider a closed subspace $V$ of $L_{\log}$ and K\"othe
function spaces $X_0,X_1\in\Cal X.$ The following lemma then
shows that Proposition 3.6 can be used.

\proclaim{Lemma 3.7}If $X\in\Cal X$ then there is a
separable K\"othe function $Z\supset X$ with $Z\in \Cal
X.$\endproclaim

\demo{Proof}Simply pick $0\le w\in X^*$ with $\log w\in L_1$
and let $Z=L_1(w).$\bull\enddemo

\vskip2truecm

\subheading{4.  Operators on K\"othe function spaces} Now
suppose that $S$ is a Polish space and that $\mu$ is a
probability measure on $S.$ We suppose that $T:L_2(\mu)\to
L_2(\mu)$ is a bounded self-adjoint operator with $\|T\|\le
1.$ Now suppose $X$ is a K\"othe function space on
$(S,\mu).$ We define $\|T\|_X=\sup\{\|Tf\|_X:\ f\in L_2\cap
X,\ \|f\|_X\le 1\}.$ If $X$ is a separable K\"othe function
space (with the Fatou property) then $L_2\cap X$ is dense in
$X$ and so $T$ extends to a bounded operator $T:X\to X$ if
and only if $\|T\|_X<\infty.$

The following remarks are elementary:

\proclaim{Lemma 4.1}(1) For any separable K\"othe function
space $X$, we have $\|T\|_X=\|T\|_{X^*}.$\newline (2) If
$X,Y$ are separable K\"othe function spaces and $0<\theta<1$
then $\|T\|_{X^{\theta}Y^{1-\theta}} \le
\|T\|_X^{\theta}\|T\|_Y^{1-\theta}.$ \endproclaim

\demo{Proof}(1) is a trivial deduction from duality and the
self-adjointness of $T$ (the Banach space adjoint of $T$ is
the complex conjugate of the Hilbert space adjoint of $T$).
(2) is immediate from complex interpolation.\bull\enddemo

Let us now say that $X$ is a {\it $T$-direction (space)} if
there exists $0<\theta<1$ so that
$\|T\|_{X^{\theta}L_2^{1-\theta}}<\infty.$ Note that if
$0<\theta<1$ then the space $X^{\theta}L_2^{1-\theta}$ is
$p$-convex and $q$-concave where $\frac1p=1-\frac1q=
\frac{1-\theta}2+\theta.$ It is thus super-reflexive and
hence separable.  Clearly, by duality, $X$ is a
$T$-direction if and only if $X^*$ is a $T$-direction.

If $w\in L_{0,\bold R}(\mu)$ we will say that $w$ is a { \it
$T$-weight direction} if there exists $\alpha>0$ so that $T$
is bounded on $L_2(e^{\alpha w}).$ Thus $w$ is $T$-weight
direction if and only if $L_2(e^w)$ is a $T$-direction.  The
space of all $T$-weight directions will be denoted by $\Cal
D=\Cal D(T).$ We define $$ \|w\|_{\Cal D}=\inf\{t>0:\
\|T\|_{L_2(e^{w/t})}\le e\}.$$ By complex interpolation it
is clear that $\|w\|_{\Cal D}<\infty$ if $w\in\Cal D.$

\proclaim{Lemma 4.2}(1) The set $\{w:\|w\|_{\Cal D}\le 1\}$
is a closed absolutely convex subset of $L_{0,\bold R}$.
\newline (2) For any $f\in L_2$ and $w\in \Cal D$ such that
$f,wf\in L_2$ we have $\|T(wf)-wTf\|_2\le e\|w\|_{\Cal
D}\|f\|_2.$ In particular $\|w\|_{\Cal D}=0$ if and only if
$T(wf)=wT(f)$ whenever $f,wf\in L_2.$ \endproclaim

\demo{Proof}Note that $\|w\|_{\Cal D}\le 1$ if and only if
$$\int |Tf|^2e^{w/t}d\mu\le e^2\int |f|^2 e^{w/t}d\mu$$
whenever $f\in L_2(1+e^w)$ and $t>1.$ It then follows from
the Dominated Convergence Theorem that $\|w\|_{\Cal D}\le 1$
if and only if $$\int |Tf|^2e^wd\mu\le e^2\int
|f|^2e^wd\mu$$ for $f\in L_2(1+e^w).$

Now suppose $w_n$ is a sequence with $w_n\to w$ a.e. and
$\|w_n\|_{\Cal D}\le 1.$ Let $u=1+\sup e^{w_n}.$ Then if
$f\in L_2(u)$ we clearly have $\|Tf\|_{L_2(e^w)}\le
e\|f\|_{L_2(e^w)}.$ By a density argument this estimate
extends to $L_2(1+e^w).$ Hence $\|w\|_{\Cal D}\le 1.$

Convexity of the set $\{w:\|w\|_{\Cal D}\le 1\}$ follows
from the fact that
$L_2(u)^{\theta}L_2(v)^{1-\theta}=L_2(u^{\theta}v^{1-\theta}).$
Symmetry follows from the fact that $L_2(u)^*=L_2(u^{-1}).$

(2) Finally suppose $\|w\|_{\Cal D}\le 1.$ Then for any real
$-1\le t\le 1$ we have $\|T\|_{L_2(e^{tw})}\le e.$ Suppose
$f,g\in \cap_{n\in\bold Z}L_2(e^{nw}).$ Then the maps $z\to
e^{zw}f$ and $z\to e^{zw}g$ are entire $L_2-$valued
functions.  It follows that the map $\varphi(z)= \int
T(e^{zw}f)e^{-zw}gd\mu$ is an entire function.  However if
$z=x+iy$ with $-1\le x\le 1,$ $$ \align |\varphi(z)| &\le
\left(\int |T(e^{zw}f|^2e^{2xw}d\mu\right)^{1/2}\left( \int
|e^{-zw}|^2 e^{2xw}|g|^2d\mu\right)^{1/2}\\ &\le
e\|f\|_2\|g\|_2 \endalign $$ and so by Cauchy's theorem,
$|\varphi'(0)|\le e\|f\|_2\|g\|_2.$ This implies that $$
|\int (T(wf)-wTf)g\,d\mu|\le e\|f||_2\|g\|_2.$$ By varying
$g$ we see that $T(wf)-wTf\in L_2$ and $\|T(wf)-wTf\|_2\le
e\|f\|_2$ whenever $f\in \cap_{n\in\bold Z}L_2(e^{nw}).$ A
simple approximation argument completes the proof that this
holds under the weaker hypothesis that $f,wf\in L_2.$

Clearly now if $\|w\|_{\Cal D}=0$ we obtain the conclusion
that $T(wf)=wTf$ under the same hypotheses.  Conversely if
$T(wf)=wTf$ for all $f$ such that $f,wf\in L_2$ it is easy
reverse the argument to show that $\|T\|_{L_2(e^{tw})}\le 1$
for all real $t.$\bull\enddemo

Now, if $X$ is a K\"othe function space, we will say that
$X$ satisfies the {\it $T$-weight condition } if there exist
constants $(C,M)$ so that if $0\le f\in B_X$ then there
exists $g\ge f$ with $\|g\|_X\le M$ and $\|\log g\|_{\Cal D}
\le C.$ We then say that $X$ satisfies the $T$-weight
condition with constants $(C,M).$

\proclaim{Theorem 4.3}(1) Suppose $X_j$ for $j=0,1$ are
K\"othe function spaces with the $T$-weight condition with
constants $(C_j,M_j).$ Then if $0<\theta<1,$
$X_0^{1-\theta}X_1^{\theta}$ has the $T$-weight condition
with constants $((1-\theta)C_0+\theta
C_1,M_0^{1-\theta}M_1^{\theta}).$ \newline (2) Suppose
$0<\theta<1$.  Then for any K\"othe function space $X,$ $X$
has the $T$-weight condition with constants $(C,M)$ if and
only if $X^{\theta}(=L_{\infty}^{1-\theta}X^{\theta})$ has
the $T$-weight condition with constants $(\theta
C,M^{\theta}).$ \endproclaim

\demo{Proof}(1)Suppose $X_j$ satisfies the $T$-weight
condition with constants $(C_j,M_j).$ Suppose $0\le f\in
X_{\theta}=X_0^{1-\theta}X_1^{\theta}$ with
$\|f\|_{X_{\theta}}\le 1.$ We may factor
$f=f_0^{1-\theta}f_1^{\theta}$ where $0\le f_j\in B_{X_j}$
for $j=1,2.$ Then pick $g_j\in X_j$ with $0\le f_j\le g_j$
and $\|g_j\|_{X_j}\le M_j$ so that $\|\log g_j\|_{\Cal D}\le
C_j.$ Then $f\le g=g_0^{1-\theta}g_1^{\theta}$ and clearly
$\|g\|_{X_{\theta}}\le M_0^{1-\theta}M_1^{\theta}$ and
$\|\log g\|_{\Cal D}\le (1-\theta)C_0+\theta C_1.$

(2) This is trivial and we omit it.  \bull\enddemo

Before proceeding we will need a technical lemma.

\proclaim{Lemma 4.4}Suppose $X$ is a K\"othe function space
with the property that there exist constants $0<c<1,\ C,M$
so that if $0\le f\in X$ there exists a Borel set $A\subset
S,$ and $g\ge f\chi_A$ such that:\newline (1)
$\|f-f\chi_A\|_X\le c\|f\|_X$ \newline (2) $\|g\|_X \le
M\|f\|_X$\newline (3) $\|\log g\|_{\Cal D}\le C.$ \newline
Then $X$ satisfies the $T$-weight condition with constants
$(C',M')$ where $C'=\max(1,C)$ and for suitable $M'.$
\endproclaim

\demo{Proof}Suppose $f=f_0\in B_X.$ We inductively define
Borel sets $(A_n)_{n=1}^{\infty},$ and sequences
$(f_n)_{n\ge 1}, \ (g_n)_{n\ge 1}$ in $X_+$ so that for
$n\ge 1,$ $$ \align \|f_{n-1}-f_{n-1}\chi_{A_{n}}\| &\le
c\|f_{n-1}\|_X,\\ g_n &\ge f_{n-1}\chi_{A_{n}}\\ \|g_n\|_X
&\le M\|f_{n-1}\|_X\\ \|\log g_n\|_{\Cal D} &\le C\\
f_{n}&=f_{n-1}-f_{n-1}\chi_{A_{n}}.  \endalign $$ Then, by
construction, $\|f_n\|_X \le c^n$ and $f-f_n=f\chi_{B_n}$
where $B_n=\cup_{k\le n}A_k.$ It follows that $f=\max_{n\ge
1}f_{n-1}\chi_{A_n}\le \max_{n\ge 1} g_n.$

Now for $0<p\le 1$ we have $$\|(\sum g_n^p)^{1/p}\|_X \le
(\sum\|g_n\|_X^p)^{1/p} \le
M(\sum_{n=0}^{\infty}c^{np})^{1/p}$$ by $p$-convexity of
$X.$ Thus $$\|(\sum g_n^p)^{1/p}\|_X \le M(1-c^p)^{-1/p}.$$
Choose $p=\min(1,1/C),$ and let $g=(\sum g_n^p)^{1/p}.$ Then
$\|\log g_n^p\|_{\Cal D}<1$ and so that if $h\in L_2(1+g)$
then $$ \int |Th|^2 g_n^pd\mu \le e^2 \int |h|^2 g_n^p
d\mu.$$ On adding we see that $$ \int |Th|^2 g^p\,d\mu \le
e^2 \int |h|^2 g^p\,d\mu.$$ Thus $\|\log g\|_{\Cal D} \le
\frac1p,$ and the result follows.\bull\enddemo

\proclaim{Theorem 4.5}Suppose $X_0$ is a K\"othe function
space satisfying the $T$-weight condition and that $X_1$ is
an arbitrary K\"othe function space.  Suppose for some
$0<\phi\le 1$ the space $X_0^{1-\phi}X_1^{\phi}$ satisfies
the $T$-weight condition.  Then for any $0<\theta<1,$ the
space $X_0^{1-\theta}X_1^{\theta}$ satisfies the $T$-weight
condition.  If $X_1$ is super-reflexive then we also have
that $X_1$ satisfies the $T$-weight condition.  \endproclaim

\demo{Proof}If $0<\theta\le \phi$ this follows immediately
from Theorem 4.3.  We therefore suppose $0<\phi<\theta<1.$
We will write $X_{\tau}=X_0^{1-\tau}X_1^{\tau}.$ Suppose
$X_0$ satisfies the $T$-weight condition with constants
$(C_0,M_0)$ and that $X_{\phi}$ satisfies the $T$-weight
condition with constants $(C_{\phi},M_{\phi}).$ We will
verify the conditions of Lemma 4.4 for the space
$X_{\theta}.$ Fix a constant $L$ so that
$L^{\frac{\phi}{\theta-\phi}}=2M_0M_{\phi}.$

Suppose $f=f_{\theta}\ge 0$ and
$\|f_{\theta}\|_{X_{\theta}}\le 1.$ Then we can write
$f_{\theta}=f_0^{1-\theta}f_1^{\theta}$ where $f_j\ge 0,$
$\|f_j\|_{X_j}\le 1$ for $j=0,1.$ Let
$f_{\phi}=f_0^{1-\phi}f_1^{\phi}$ so that
$\|f_{\phi}\|_{X_{\phi}}\le 1.$ Then there exists
$g_{\phi}\ge f_{\phi}$ with $\|g_{\phi}\|_{X_{\phi}}\le
M_{\phi}$ and $\|\log g_{\phi}\|_{\Cal D}\le C_{\phi}.$

We thus write $g_{\phi}=g_0^{1-\phi}g_1^{\phi}$ where
$g_j\ge 0$ and $\|g_j\|_{X_j}\le M_{\phi}.$ Then there
exists $h_0\ge g_0$ with $\|h_0\|\le M_0M_{\phi}$ and
$\|\log h_0\|_{\Cal D} \le C_0.$ Next we define
$h_{\theta}=h_0(g_{\phi}h_0^{-1})^{\alpha}$ where
$\alpha=\theta/\phi>1$ and $(0/0)=0.$ Then
$h_{\theta}=g_{\phi}(g_{\phi}h_0^{-1})^{\alpha-1} \le
g_0^{1-\theta}g_1^{\theta}.$ Thus
$\|h_{\theta}\|_{X_{\theta}}\le M_{\phi}.$ Let
$h'=Lh_{\theta}.$ We note that $\|h'\|_{X_{\theta}}\le
LM_{\phi}$ and $$ \align \|\log h'\|_{\Cal D}&=\|\log
h_{\theta}\|_{\Cal D}\\ &= \|(1-\alpha)\log h_0 +\alpha\log
g_{\phi}\|_{\Cal D}\\ &\le (\alpha-1)C_0 + \alpha C_{\phi}.
\endalign $$

Let $A=\{s:\ f(s)\le Lh(s)\}$ and let $B=S\setminus A.$ Then
$f\chi_A\le h'$ and $\|f-f\chi_A\|_{X_{\theta}} \le
\|f_0\chi_B\|_{X_0}^{1-\theta}.$ Now if $s\in B$ we have $$
\align \frac{f(s)}{f_{\phi}(s)} &\ge
L\frac{h_{\theta}(s)}{g_{\phi}(s)}\\ &=
L\left(\frac{g_{\phi}(s)}{h_0(s)}\right)^{\alpha-1}.
\endalign $$ Hence $$ \align f_0(s) &=
f_{\phi}(s)\left(\frac{f_{\phi}(s)}{f(s)}\right)^{\frac{-\phi}{\theta-\phi}}\\
&\le L^{-\frac{\phi}{\theta-\phi}}g_{\phi}(s)\left(
\frac{g_{\phi}(s)}{h_0(s)}\right)^{-1}\\ &\le \frac12
M_0^{-1}M_{\phi}^{-1}h_0(s).  \endalign $$ We conclude that
$\|f_0\chi_B\|_{X_0}\le \frac12$ and hence that
$\|f-f\chi_A\|_{X_{\theta}} \le (\frac12)^{1-\theta}=c<1$
say.  Thus the hypotheses of Lemma 4.4 are verified and
$X_{\theta}$ has the $T$-weight condition.

For the last assertion, if $X_1$ is super-reflexive, we may
suppose that there is a K\"othe function space $Y$ so that
$X_1=X_0^{1-\tau}Y^{\tau}$ for $0<\tau<1$.  The above
argument then gives the conclusion.  \bull\enddemo

Let us draw a simple conclusion form Theorem 4.5.

\proclaim{Theorem 4.6}Suppose $X$ is $q$-concave for some
$q<\infty$ and that $\Phi_X=\sum_{j=1}^n\alpha_j\Phi_{X_j}$
where $\alpha_j\in\bold R$ and each $X_j$ satisfies the
$T$-weight condition.  Then $X$ satisfies the $T$-weight
condition.\endproclaim

\demo{Proof}Since we may replace $X$ by $X^{\alpha}$ where
$0<\alpha<1$ we consider only the case when $X$ is
super-reflexive.  It clearly also suffices to establish this
theorem when $n=2.$ It follows directly from Theorem 4.3
when $\alpha_1,\alpha_2\ge 0.$ If $\alpha_1,\alpha_2\le 0$
then $-\Phi_X$ is convex so that $\Phi_X$ is linear when
$X=wL_{\infty}$ for some weight $w$ which contradicts
super-reflexivity.  We may thus suppose $\alpha_1$ and
$\alpha_2$ have opposite signs and by Theorem 4.3 we need
only consider the case $\alpha_1=1$ and $\alpha_2<0.$ Define
$Y_1=X^{1/2}$ and then let $Y_0$ be defined by
$Y_0=X_1^{1/2}X_2^{1/4}.$ Then $Y_0$ is an interpolation
space between $X$ and $X_1^{1/2}X_2^{1/2}$ and so is
super-reflexive.  By Theorem 4.3, $Y_0$ satisfies the
$T$-weight condition; but for an appropriate $\phi>0$ we
have $Y_0^{1-\phi}Y_1^{\phi}=X_1^{1/2}$ which also satisfies
the $T$-weight condition.  Now by Theorem 4.6 we conclude
that $X^{1/2}$ satisfies the $T$-weight condition and thus
Theorem 4.3 completes the proof.\bull\enddemo

\proclaim{Lemma 4.7}Let $X$ be an exactly 2-convex K\"othe
function space and let $Y=(X^2)^*.$ Then:\newline (1) If $Y$
satisfies the $T$-weight condition with constants $(1,M)$
then $\|T\|_X<\infty.$ \newline (2) If $\|T\|_X<\infty$ then
$Y$ satisfies the $T$-weight condition.  \endproclaim

\demo{Proof}(1) Suppose $f\in L_2\cap B_X$.  Suppose $0\le
u\in Y$ with $\|u\|_Y\le 1.$ Then there exists $v\ge u$ so
that $\|v\|_Y\le M$ and $\|\log v\|\le 1.$ Thus (cf.  [30],
Theorem A'), $$ \int |Tf|^2 u\, d\mu \le e^2 \int
|f|^2v\,d\mu\le Me^2\||f|^2\|_{X^2} $$ and hence $\|Tf\|_X
\le M^{1/2}e\|f\|_X.$

(2) This follows from a result of Rubio de Francia ([30],
Theorem A'); in the case when $X$ is a weighted $L_p-$space
for $p>2$ it was shown by Cotlar and Sadosky [9].  In fact
by Theorem A' of [30] there is a constant $M$ so that if
$0\le u\in B_Y,$ there exists $v\ge u$ with
$\|T\|_{L_2(v)}\le M$ and $\|v\|_Y\le 2\|u\|_Y.$ Now by
interpolation $\|\log v\|_{\Cal D}\le \log M$ and we are
done.\bull\enddemo

\proclaim{Lemma 4.8}Let $X$ be an exactly 2-convex K\"othe
function space and let $Y=(X^2)^*.$ Then $Y$ satisfies the
$T$-weight condition if and only if $X$ is a $T$-direction
space.  \endproclaim

\demo{Proof}Suppose $Y$ satisfies the $T$-weight condition
with constants $(C,M).$ Choose $\theta>0$ so that $\theta
C<1.$ Consider the space $Z=L_2^{1-\theta}X^{\theta}$.  Then
$\Phi_Z=\frac{(1-\theta)}2\Lambda+\theta\Phi_X$ and so
$2\Phi_Z +\theta\Phi_Y=\Lambda$ so that
$Y^{\theta}=(Z^2)^*.$ By applying Lemma 4.7, $T$ is bounded
on $Z.$

Conversely if $X$ is a $T$-direction space there exists
$\theta>0$ so that $\|T\|_Z<\infty$ where
$Z=L_2^{1-\theta}X^{\theta}$ and so by Lemma 4.7,
$Y^{\theta}$ satisfies the $T$-weight condition.  Theorem
4.3 completes the proof.\bull\enddemo

We are now finally able to state our main result of this
section.

\proclaim{Theorem 4.9}Suppose $T:L_2\to L_2$ is a
self-adjoint operator with $\|T\|\le 1.$ Suppose
$L_{\infty}$ is a $T$-direction space (i.e. there exists
$p>2$ so that $\|T\|_{L_p}<\infty$).  Then \newline (1) If
$X$ is satisfies the $T$-weight condition then $X$ is a
$T$-direction space.\newline (2) If $X$ has is $q$-concave
for some $q<\infty$ then $X$ is a $T$-direction space if and
only if $X$ satisfies the $T$-weight condition.
\endproclaim

\demo{Remark}Note that $L_{\infty}$ is always a $T$-weight
space.  In general our assumption that $T$ is bounded at
some $L_p$ where $p>2$ is equivalent to the requirement that
$L_2$ satisfies the $T$-weight condition by Lemma 4.8.  This
shows that the assumption is necessary for the theorem to
hold.\enddemo

\demo{Proof}We assume that $p>2$ and $p'<2$ are conjugate
indices so that $\|T\|_{L_p}=\|T\|_{L_{p'}}<\infty.$ We
first notice that $L_2$ must satisfy the $T$-weight
condition.  Indeed by Lemma 4.7, $L_r$ satisfies the
$T$-weight condition when $\frac1r+\frac2p=1$ and hence by
Theorem 4.3, $L_2$ satisfies the $T$-weight condition.

We will now prove (2) under the stronger hypothesis that $X$
is super-reflexive.

We next show that, in general, if $X$ is super-reflexive and
satisfies the $T$-weight condition, then $X^*$ also
satisfies the $T$-weight condition.  In fact
$L_2=X^{1/2}(X^*)^{1/2}$ and so it follows from Theorem 4.5
that $X^*$ has the $T$-weight condition.

We now proceed to the proof of the theorem.  Assume first
$X$ is super-reflexive and satisfies the $T$-weight
condition.  We now may select $0<\alpha<1$ small enough so
that $(X^*)^{\alpha}$ has the $T$-weight condition with
constants $(1,M)$ for suitable $M.$ Now by Lemma 4.7 $T$ is
bounded on the space $Z$ where
$2\Phi_Z=\Lambda-\alpha\Phi_{X^*},$ or
$\Phi_Z=\frac12(1-\alpha)\Lambda+\frac12\alpha\Phi_X.$ Now
$\Delta_{\Phi_Z}\ge \frac12(1-\alpha)\Delta_{\Lambda}$ so
that $Z$ has nontrivial concavity and is thus
super-reflexive.  It follows that $T$ is also bounded on any
space $Y=L_{p'}^{1-\beta}Z^{\beta}$ where $0<\beta<1.$ We
select $\beta$ so that $Y$ is an interpolation space between
$L_2$ and $X.$ In fact $$ \Phi_Y =((1-\beta) (1-\frac1p)
+\frac{\beta}{2}(1-\alpha))\Lambda
+\frac12\alpha\beta\Phi_X$$ and the conclusion is obtained
by choosing $\beta$ so that $$
2(1-\beta)(1-\frac1p)+{\beta}(1-\alpha) +\frac12 \alpha\beta
=1$$ or $$ \beta=\frac{1-\frac2p}{\frac\alpha2+1-\frac2p}.$$
Now $T$ is bounded at $Y$ and so $X$ is a $T$-direction
space.

Now suppose, conversely that $X$ is a $T$-direction space.
Then for suitable $\theta>0,$ $T$ is bounded at
$L_2^{1-\theta}X^{\theta}.$ Interpolating with $L_p$ we see
that $T$ is also bounded at any space $Z$ where
$$\Phi_Z=(\frac{1-\alpha}p +\frac{\alpha}2(1-\theta))\Lambda
+\theta\alpha\Phi_X$$ where $0<\alpha<1.$ Notice that
$$\Delta_{\Phi_Z} \le (\frac{1-\alpha}p
+\frac{\alpha(1+\theta)}2)\Delta_{\Lambda}$$ and so by
choosing $\alpha$ small enough we can suppose that $Z$ is
2-convex.  Let us put $\Phi_Z=\beta\Lambda+\gamma\Phi_X$
where $0<\beta,\gamma$ and $\beta +\gamma<\frac12.$ Thus $Y$
satisfies the $T$-weight condition where $\Phi_Y
=\Lambda-2\Phi_Z.$ We can now solve for $\Phi_X$ in form
$$\Phi_X=\frac{1}{2\gamma}((1-2\beta)\Lambda -\Phi_Z).$$ An
application of Theorem 4.6 now completes the proof for the
case when $X$ is super-reflexive.

Now consider (1).  If $X$ satisfies the $T$-weight
condition, then so does $L_2^{1/2}X^{1/2}$ by Theorem 4.6
(or 4.5).  This space is super-reflexive and so is also a
$T$-direction space; hence $X$ is $T$-direction space.

Finally we complete the proof of (2) when $X$ is $q$-concave
for some finite $q$ and is a $T$-direction space.  Then
$Y=L_2^{1/2}X^{1/2}$ is a $T$-direction space and is
super-reflexive; hence it satisfies the $T$-weight
condition.  Since $\Phi_X= 2\Phi_Y-\frac12\Lambda$, we
complete the proof by Theorem 4.6.\bull\enddemo

\vskip2truecm

\subheading{5.  Interpolation of Hardy spaces} We again
consider a probability measure $\mu$ on a Polish space $S.$
Consider the Orlicz algebra $L_{\log}$, and let $\Cal X$ be
the collection of all K\"othe functions spaces $X$ so that
$X,X^*\in \Cal X.$ Consider a closed subalgebra $H$ of
$L_{\log}$ (which is always assumed to contain the
constants).  We define for every $X\in\Cal X$ the Hardy
space $H_X=H\cap X$ so that $H_X$ is a closed subspace of
$X.$ In particular we define $H_p=L_p\cap X$ when $1\le p\le
\infty.$

We will say that $H$ is of {\it Dirichlet type} if for every
invertible $f\in L_{\log}$ there exists $g\in H$ which is
invertible in $H$ so that $|g|=|f|$ a.e.; equivalently, $H$
is a Dirichlet-type algebra if for every real $u\in L_1$
there exists an invertible $g\in H$ with $|g|=e^u$ a.e.

The simplest example of such a Dirichlet-type algebra is the
Smirnov class $N^+$ (or Hardy algebra) considered as a
subalgebra of $L_{\log}(\bold T).$ In this way one generates
the standard Hardy spaces.  More generally suppose $A$ is a
subalgebra of $L_{\infty}(S,\mu)$ so that that $f\to \int
f\,d\mu$ is a multiplicative linear functional and $\Re A$
is weak$^*$-dense in $L_{\infty,\bold R}.$ Thus $A$ is a
{\it weak$^*$-Dirichlet algebra} (cf.  [1], [12], [14]).
Let $H$ be the closure of $A$ in $L_{\log}$; then $H$ has
the Dirichlet property and the standard abstract Hardy
spaces are obtained.  The reader may consult Gamelin [12]
for details when $A$ is generated by a Dirichlet algebra;
see also Barbey-K\"onig [1].

Another example is obtained when one consider $(\bold
T\times S,\lambda\times \mu)$ and defines $H$ to be the
space of all functions $f(t,s)$ so that $f\in L_{\log}$ and
for a.e.  $s\in S$ the function $f_s\in N^+$ where
$f_s(t)=f(t,s).$ In this way we can treat vector-valued
problems.

Notice that, in each case, one can always replace the
measure $\mu$ by a measure $w\,d\mu$ as long as $w,\log w\in
L_1$.  This will not change $L_{\log}$ or $H$ but will alter
the space $H_2$.  This change of density allows one to study
skew projections.

\proclaim{Lemma 5.1}Let $H$ be any closed subalgebra of
$L_{\log}.$ If $f\in H_1$ then $e^f\in H.$\endproclaim

\demo{Proof}The series $\sum_{n\ge 0}\frac{f^n}{n!}$
converges in $L_{\log}$ since it converges a.e. and
$\sum_{n\ge 0}\frac{|f|^n}{n!}=e^{|f|}\in
L_{\log}.$\bull\enddemo

\proclaim{Lemma 5.2}Suppose $H$ is a Dirichlet-type algebra.
Then if $f\in H$ and $v\in L_1$ there is a sequence $g_n\in
H$ so that $|g_n|\le \min(ne^v,|f|)$ and $g_n\to f$ in
measure (and hence in $L_{\log}).$ \endproclaim

\demo{Proof}First consider the subspace $G$ of $L_{1,\bold
R}\times L_{0,\bold R}$ of all $(u,v)$ so that
$e^{\pm(u+iv)}\in H.$ It is easy to check that $G$ is
closed.  Hence by an application of the Open Mapping Theorem
if $\|u_n\|_1\to 0$ there exist $v_n\to 0$ in $L_0$ so that
$e^{u_n+iv_n}\in H.$

Now pick any $h\in L_{1,\bold R}$ with $h\ge |f|.$ There
exists an invertible $\psi\in H$ with $|\psi|=e^h.$ Now let
$u_n=h-\min(h,v+\log n).$ Then $\|u_n\|_1\to 0$ and so there
exist $v_n\to 0$ in measure so that $e^{\pm(u_n+iv_n)}\in
H.$ Let $g_n=\psi e^{-(u_n+iv_n)}f$ and the result follows
easily.\bull\enddemo

Suppose $H$ is a closed subalgebra of $L_{\log}.$ We define
$V$ to be the subspace of $L_{\log}$ of all $f$ so that
$\int \bar fg\,d\mu=0$ whenever $g\in H$ and $fg\in L_1.$
For $X\in\Cal X$ we set $V_X=V\cap X.$ It is trivial to see
that if $f\in V$ and $g\in H$ then $f\bar g\in V.$

\proclaim{Lemma 5.3}Assume $H$ is a Dirichlet-type
algebra.\newline (1) $f\in V$ if and only if there exists an
invertible $g\in H$ so that $fg\in L_1$ and $\int \bar
fgh\,d\mu=0$ for every $f\in H_{\infty}.$\newline (2) $V$ is
a closed subspace of $L_{\log}.$ \newline (3) If $X\in\Cal
X$ then $X\cap H_{\infty}$ is dense in $H_X$ and $X\cap
V_{\infty}$ is dense in $V_X.$ \endproclaim

\demo{Proof} (1) Suppose $\psi\in H$ and $f\psi\in L_1.$
Then by Lemma 5.2 there exist $\psi_n\in H$ so that
$|\psi_n|\le \min(n|g|,|\psi|)$ and $\psi_n\to\psi$ in
measure.  Then $\int \bar f \psi_n d\mu=0$ and the
conclusion follows from Dominated Convergence.

(2) Suppose $f_n\to f$ in $L_{\log}$ where $f_n\in V.$ By
passing to a subsequence we can suppose that
$F=\sup_n|f_n|\in L_{\log}.$ Choose any invertible $g\in H$
so that $|g|\ge F$.  Then $\int \bar fg^{-1}h\,d\mu=0$ for
every $h\in H_{\infty}.$ Hence by (1), $f\in V.$

(3) Suppose $f\in X$; then (Lemma 2.2) there exists $w\ge
|f|$ with $\log w\in L_1.$ Then there exists an invertible
$g\in H$ with $|g|=w$ a.e. and by Lemma 5.2 a sequence
$g_n\in H_{\infty}$ with $|g_n|\le |g|$ so that $g_n\to g$
in measure.  If $f\in H$ then the sequence $(fg^{-1}g_n)$ is
in $H_{\infty}\cap X$, converges in measure to $f$ and is
lattice bounded by $|f|.$ Hence it converges also in $X.$ If
$f\in V$ we use a similar argument on $f\bar g^{-1}\bar
g_n.$\bull\enddemo

>From now on, we suppose $H$ is a Dirichlet-type algebra.  We
define $\Cal R$ to be the orthogonal projection of $L_2$
onto $H_2$; it follows from the preceding lemma that the
kernel of $\Cal R$ is $V_2.$ Further, if $X\in\Cal X,$ then
$\Cal R$ is bounded at $X$ if and only if $X=H_X\oplus V_X.$
We will say that $H$ is a {\it Hardy-type algebra} if
$L_p=H_p\oplus V_p$ for all $1<p<\infty.$ Note that all the
examples quoted are of Hardy type.

If $w\in L_{1,\bold R}$ we will say that $w\in BMO$ if $w\in
H_1+L_{\infty}$ and we define the BMO-norm by
$\|w\|_{BMO}=\inf\{\|w-h\|_{\infty}:h\in H_1\}.$ Let us note
in passing that the infimum is attained.  Indeed if $h_n\in
H_1$ is such that $\|w-h_n\|_{\infty}\to \|w\|_{BMO}$ then
by Komlos's theorem [22], since $(h_n)$ is $L_1$-bounded, we
can pass to a sequence of convex combinations $(g_n)$ of
$(h_n)$ which converge a.e. to some $g.$ However it is
easily seen that $\|g_n-g\|_p\to 0$ when $p<1$ and so $g\in
H$ since $H$ is closed in $L_{\log}.$

\proclaim{Proposition 5.4} If $w\in L_1,$ then $w\in BMO$ if
and only if $w$ is an $\Cal R-$weight direction.  Further
there is a constant $C$ so that if $w\in BMO$ then
$C^{-1}\|w\|_{BMO}\le \|w\|_{\Cal D(\Cal R)}\le
C\|w\|_{BMO}.$\endproclaim

\demo{Proof}First suppose $w\in\Cal D\cap L_1$.  By Lemma
4.2, if $f\in V_2$ then $\|\Cal R(wf)\|_2\le e\|w\|_{\Cal
D}\|f\|_2.$ Now suppose $f\in V_{\infty},$ with $wf\in L_2.$
Then for $\epsilon>0$ there exists an invertible $g\in H$ so
that $|g|=|f|^{1/2}+\epsilon$ a.e.  Then $f\bar g^{-1}\in
V_2$ and $g\in H_2$.  $$ \int fw\,d\mu = \int \bar g(wf\bar
g^{-1})d\mu =\int \bar g\Cal R(wf\bar g^{-1})d\mu.$$ Hence
we have $$ |\int fw\,d\mu| \le e\|w\|_{\Cal
D}\|fg^{-1}\|_2\|g\|_2.$$ Letting $\epsilon\to 0$ we obtain
$$ |\int fw\,d\mu| \le e\|w\|_{\Cal D}\|f\|_1.$$ Now by the
Hahn-Banach theorem there exists $\psi\in L_{\infty}$ so
that $\|\psi\|_{\infty}\le e\|w\|_{\Cal D}$ and $\int
f(w-\bar\psi)d\mu=0$ for $f\in V_{\infty}\cap L_2(|w|^{2}).$
Now for any $f\in V_{\infty}$ we can find, utilizing Lemma
5.2, with $v=-(\log_+|f|+\log_+|w|)$ an invertible $g_n\in
H$ so that $|g_n|\le \min(1,n|f|^{-1}w^{-1})$ and $g_n\to 1$
a.e.  Then $\bar g_nf\in V_{\infty}\cap L_2(w^2)$ and by the
Dominated Covergence Theorem we have $\int
f(w-\bar\psi)d\mu=0.$

It now follows, again from the Hahn-Banach theorem that
$w-\psi\in H_1$ and hence that $\|w\|_{BMO}\le e\|w\|_{\Cal
D}.$

Now conversely suppose $w \in BMO,$ with $\|w\|_{BMO}\le 1.$
Let $X_0=L_2(e^{2w})$ and $X_1=L_2(e^{-2w}).$ Then if
$X_{\theta}=[X_0,X_1]_{\theta}$ we have $X_{1/2}=L_2.$ We
claim that $H$ is interpolation stable at $1/2$ and further
there is a universal constant $C$ so that $K(\frac12,H)\le
C.$ In fact there exists $h\in H_{\infty}$ so that
$\|w-h\|_{\infty}\le 1.$ Suppose $f\in H_2.$ Then we define
a map $F:\Cal S \to H$ by $$F(z)=e^{-1+4z^2} e^{(1-2z)h}f.$$
It is clear that $F$ is analytic into $H$ and $$ \align \int
|F(it)|^2e^{-2w}d\mu &=e^{-2-8t^2}\int |e^{2(1+2it)(h-w)}|
|f|^2d\mu\\ &= e^{4|t|-8t^2} \|f\|_2^2\le e^{1/2}\|f\|_2^2
\endalign $$ while $$ \align \int |F(1+it)|^2e^{-2w}d\mu
&\le e^{6-8t^2}\int |e^{2(-1-2it)(h-w)}||f|^2d\mu\\ &\le
e^{4+4|t|-8t^2}\|f\|_2^2 \le e^{9/2}\|f\|_2^2.  \endalign $$
It follows that $H$ is interpolation stable at $1/2$ with
$K(1/2,H)\le e^{5/4}.$ Now it follows from Theorem 3.3 and
its proof that $L_2(e^{w/t})=H_2(e^{w/t})\oplus
V_2(e^{w/t})$ if $|t|\le C$ for some absolute constant $C$.
Thus $w\in\Cal D(\Cal R)$ and $\|w\|_{\Cal D}\le
C.$\bull\enddemo

We will now say that a K\"othe function space $X\in \Cal X$
is {\it BMO-regular} (for $H$) if there are constants
$(C,M)$ so that if $0\le f\in X$ with $\|f\|_X\le 1$ then
there exists $g\in X$ with $g\ge f,$ $\|g\|_X\le M$ and
$\|\log g\|_{BMO}\le C.$

\proclaim{Lemma 5.5}Suppose $X\in\Cal X$.  Then $X$
satisfies the $\Cal R-$weight condition if and only if $X$
is BMO-regular.  \endproclaim

\demo{Proof}One direction is obvious.  For the other, note
that if $X$ the $\Cal R-$weight direction, then given $f\in
X_+$ with $\|f\|_X=1$ there exists $f'\ge f$ with
$\|f'\|_X\le 2$ and $\log f'\in L_1$ by Lemma 2.2.  Thus if
$X$ satisfies the $\Cal R$-weight condition with constants
$(C,M)$ there exists $g\ge f'$ with $\|g\|_X\le 2M$ and
$\|\log g\|_{\Cal D}\le C.$ But then also $\log g\in L_1$ so
that $\|\log g\|_{BMO}\le C'$ for a suitable constant
$C'.$\enddemo

\proclaim{Proposition 5.6}Let $H$ be a Hardy-type algebra.
If $X\in\Cal X$ is super-reflexive then $X$ is BMO-regular
if and only if $X$ is an $\Cal R-$direction.  In particular
each $L_p$ is BMO-regular.  \endproclaim

\demo{Proof}This is simply Theorem 4.9.  \bull\enddemo

\proclaim{Theorem 5.7}Suppose $H$ is a Dirichlet-type
algebra.  Suppose $X_0,X_1\in \Cal X$ are both BMO-regular
and that $0<\theta<1.$ Let
$X_{\theta}=X_0^{1-\theta}X_1^{\theta}.$ Suppose either that
(a) both $X_0,X_1$ are separable or (b) $X_{\theta}$ is
reflexive.  Then $H$ is interpolation stable at $\theta$ for
$(X_0,X_1)$ i.e.
$[H_{X_0},H_{X_1}]_{\theta}=H_{X_{\theta}}$ \endproclaim

\demo{Proof}We suppose that for $j=0,1,$ $X_j$ are
BMO-regular with constants $(C_j,M_j).$ Suppose $f\in
H_{X_{\theta}}$ with $\|f\|_{X_{\theta}}=1$; then we can
factor $|f|=f_0^{1-\theta}f_1^{\theta}$ where $0\le f_0,f_1$
and $\|f_j\|_{X_j}=1$ for $j=0,1.$ Pick $f'_j\ge f_j$ so
that $\|f'_j\|_{X_j}\le M_j$ and $\|\log f'_j\|_{BMO}\le
C_j.$ Then pick $h_j\in H_1$ so that $\|\log
f'_j-h_j\|_{\infty}\le C_j.$ We consider the following
function for $z\in\Cal S,$ $$ F(z) =
e^{z^2-\theta^2}e^{(z-\theta)(h_1-h_0)}f.$$ $F$ is
continuous into $H$ and $F(\theta)=f.$ Further if $z=j+it$
where $j=0,1$ $$ |F(j+it)| \le
|e^{j^2-\theta^2-t^2}|e^{(|j-\theta|+|t|)(C_0+C_1)}|f'_j|.$$
Hence $F\in \Cal F(X_0,X_1;H)$ (see Section 3).  Thus we get
an estimate $$ \|F(j+it)\|_{X_j}\le C'$$ where
$C'=C'(C_0,C_1,M_0,M_1,\theta).$ We can now appeal to
Proposition 3.6 to deduce that
$\|f\|_{[H_{X_0},H_{X_1}]_{\theta}}\le C'.$ and this proves
the theorem.  \enddemo

\demo{Remark} In the case when $S=\bold T$ and $H=N^+$ then
the spaces $L_p$ satisfy the BMO-condition.  This immediate
from Proposition 5.6 but there is an amusing alternative
argument.  It suffices to consider the case $p=2.$ The
Hardy-Littlewood maximal function $\Cal M$ is bounded on
$L_2$ (cf.  [31]) and for any $f\in L_2,$ $\log \Cal Mf\in
BMO$ by a result of Coifman-Rochberg [6] with an appropriate
bound.  Combining these facts shows that $L_2$ and every
$L_p$ satisfies the BMO-condition.  Notice that this then
implies an immediate proof of a well-known theorem of P.
Jones [16], [17] that $[H_{\infty},H_1]_{\theta}=H_p$ where
$p=\frac1\theta.$ Interpolation with $H_p$ when $p<1$ can be
handled in the same way.

\enddemo

To understand the picture for interpolation in general, we
need two further lemmas:

\proclaim{Lemma 5.8}Suppose $X_0,X_1\in \Cal X$ are
separable K\"othe function spaces and $0<\theta<1$ is such
that $H$ is interpolation stable at $\theta$ for
$(X_0,X_1).$ Suppose $Y_0,Y_1\in\Cal X$ are also separable
K\"othe function spaces so that for K\"othe function space
$W$ we have $Y_j=X_jW$ for $j=0,1.$ Then $H$ is
interpolation stable for $(Y_0,Y_1).$\endproclaim

\demo{Proof}As usual let
$X_{\theta}=X_0^{1-\theta}X_1^{\theta}$ and
$Y_{\theta}=Y_0^{1-\theta}Y_1^{\theta}.$ Suppose $K$ is the
constant of interpolation stability at $\theta$ for
$(X_0,X_1).$ Suppose $f\in H_{Y_{\theta}},$ and
$\|f\|_{Y_{\theta}}=1.$ Then we can factorize $f=bw$ where
$\|b\|_{X_{\theta}}=\|w\|_{W}=1.$ Now pick $b'\ge |b|$ so
that $\log b'\in L_1$ and $\|b'\|_{X_{\theta}}\le 2.$ Then
there exists an invertible $g\in H$ with $|g|=|b'|.$ Hence
there exists an admissible $F:\Cal S\to H$ so that
$F(\theta)=g$ and (a.e.)  $\|F(j+it)\|_{X_j}\le 2K.$ Define
$G:\Cal S\to H$ by $G(z)=F(z)fg^{-1}.$ It is easy to see
that $\|G(j+it)\|_{Y_j}\le 2K$ (a.e.) and $G(\theta)=f.$ $H$
is stable at $(Y_0,Y_1)$ for $\theta.$ \enddemo

\proclaim{Lemma 5.9}Suppose $X_0,X_1\in\Cal X$ are separable
K\"othe function spaces such that $H$ is interpolation
stable at $\theta$ for $(X_0,X_1).$ Then $V$ is also
interpolation stable at $\theta.$\endproclaim

\demo{Proof}Suppose $f\in
X_{\theta}=X_0^{1-\theta}X_1^{\theta}$ and $f\in V$ with
$\|f\|_{X_{\theta}}=1.$ Pick any $f'\in X_{\theta}$ so that
$f'\ge |f|,$ $\|f'\|_{X_{\theta}}\le 2$ and $\log f'\in
L_1.$ Then pick $g\in H$ so tha $|g|=f'.$ There exists an
admissible $F:\Cal S \to H$ so that $F(\theta)=g$ and
$\|F(j+it)\|_{X_j}\le 2K$ almost everywhere.  Define $\bar
F(z)=\overline{F(\bar z)}$ and consider $G(z)=f\bar
g^{-1}\bar F(z).$ Then $G$ is also admissible but has range
in $V$, $G(\theta)=f$ and $\|G(j+it)\|_{X_j}\le 2K$ a.e. so
$V$ is also interpolation stable at $\theta.$\bull\enddemo

\proclaim{Lemma 5.10}Suppose $X_0,X_1\in\Cal X$ are
separable K\"othe function spaces so that $H$ is
interpolation stable at some $0<\theta<1$ for $(X_0,X_1).$
If $\Cal R$ is bounded at
$X_{\theta}(=X_0^{1-\theta}X_1^{\theta})$ then there exists
$\eta>0$ so that $\Cal R$ is also bounded on $X_{\phi}$ if
$|\phi-\theta|\le \eta.$ \endproclaim

\demo{Proof}This follows directly from Theorem 3.3 and Lemma
5.9.\bull\enddemo

\demo{Remark}Let us note that this implies that if $L_2$ is
BMO-regular then since $H$ must be interpolation-stable at
$\theta=\frac12$ for $(L_{3/2},L_3)$ then $\Cal R$ is
bounded on $L_p$ for some $p>2.$ This provides a weak
converse to Proposition 5.6.

For the remainder of this section we require that $H$ is of
Hardy type i.e. the Riesz projection is bounded on $L_p$ for
$1<p<\infty.$

\proclaim{Proposition 5.11}Suppose $H$ is of Hardy type and
$X\in\Cal X$ is $q$-concave for some $q<\infty.$ Then $X$ is
a $\Cal R$-direction space if and only if $X$ is
BMO-regular.  \endproclaim

\demo{Proof}By Theorem 4.9 and Lemma 5.5 we obtain the
result for super-reflexive $X.$ In the general case if $X$
is BMO-regular then so is $L_2^{1/2}X^{1/2}$ and this must
therefore be a $\Cal R$-direction space, which implies $X$
is an $\Cal R-$direction space.  Conversely if $X$ is an
$\Cal R$-direction space then $L_2^{1/2}X^{1/2}$ is
BMO-regular.  But then Theorem 4.6 implies $X^{1/2}$ is
BMO-regular since it is super-reflexive.  This in turn
implies $X$ is BMO-regular.\bull\enddemo

To state our main theorem we introduce the idea of a
BMO-direction.  If $X_0,X_1\in\Cal X$ we define a K\"othe
function space $Z$ by
$\Phi_Z=\frac12(\Lambda+\Phi_{X_1}-\Phi_{X_0}).$ We say that
$X_0\to X_1$ is a {\it BMO-direction} if $Z$ is an $\Cal
R$-direction space.  If either $X_0$ is $p$-convex where
$p>1$ or $X_1$ is $q$-concave where $q<\infty$ then $Z$ has
nontrivial concavity and so this is the same as requiring
that $Z$ is BMO-regular.  If, for example, both spaces are
super-reflexive, and $X_0$ is already BMO-regular then
$X_0\to X_1$ is a BMO-direction if and only if $X_1$ is
BMO-regular; this follows immediately from Theorem 4.6.  On
an intuitive level, $X_0\to X_1$ is a BMO-direction if and
only if the parallel complex interpolation scale through
$L_2$ only passes through BMO-regular spaces.

\proclaim{Theorem 5.12}Suppose $H$ is a Hardy-type algebra
and that $X_0,X_1\in\Cal X$ are super-reflexive K\"othe
function spaces.  Then for any $0<\theta<1,$ $H$ is
interpolation stable at $\theta$ for $(X_0,X_1)$ if and only
if $X_0\to X_1$ is a BMO-direction.

In particular, if $X_0$ is BMO-regular then $H$ is
interpolation stable at $\theta$ for $(X_0,X_1)$ if and only
if $X_1$ is BMO-regular.\endproclaim

\demo{Proof}We may suppose that both $X_0,X_1$ are
$p$-convex and $q$-concave (with constant one) where
$\frac1p+\frac1q=1,$ and $1<p\le q<\infty.$ Let
$\epsilon=\frac1{2q}.$ Let $X_{\tau}=X_0^{1-\tau}X_1^{\tau}$
for $0<\tau<1.$

We start with some remarks on the implications of $H$ being
of Hardy type.  In this situation we can apply Theorem 4.9
to $\Cal R$ to deduce that a super-reflexive $X\in\Cal X$ is
BMO-regular if and only if $X$ is a $\Cal R-$direction
space.  Note that $L_2$ is BMO-regular and further that $X$
is BMO-regular if and only if $X^*$ is BMO-regular,

Let us first suppose that $H$ is interpolation stable at
$\theta.$ Now $X_{\theta}$ is $p$-convex; furthermore there
is a K\"othe function space $W$ defined by
$\Phi_W=\frac1p\Lambda-\Phi_{X_{\theta}}.$ Now consider the
g-convex K\"othe spaces $Y_{\phi}$ defined by
$\Phi_{Y_{\phi}}=\frac1p\Lambda+\Phi_{X_{\phi}}-\Phi_{X_{\theta}}.$
Clearly $$|\Delta_{\Phi_{Y_{\phi}}}-
\frac1p\Delta_{\Lambda}|\le|\phi-\theta|\Delta_{\Lambda}.$$
Hence if $|\phi-\theta|\le \epsilon=\frac1{2q}$ then
$Y_{\phi}$ is a super-reflexive K\"othe function space.  We
set $\phi_0=\theta-\epsilon$ and $\phi_1=\theta+\epsilon.$
Then $L_p=Y_{\theta}=Y_{\phi_0}^{1/2}Y_{\phi_1}^{1/2}.$

Now $H$ is interpolation stable at $\frac12$ for
$(X_{\phi_0},X_{\phi_1})$ since it is also interpolation
stable at $\theta$ for $(X_0,X_1)$.  By Lemma 5.8 $H$ is
interpolation stable at $\frac12$ also for
$(Y_{\phi_0},Y_{\phi_1}).$ However
$Y_{\phi_0}^{1/2}Y_{\phi_1}^{1/2}=L_p$ and $\Cal R$ is by
assumption bounded at $L_p.$

By Lemma 5.10, we conclude that there exists $\eta>0$ so
that $Y_{\phi}=H_{Y_{\phi}}\oplus V_{Y_{\phi}}$ for
$|\phi-\theta|\le2\eta.$ In particular $\Cal R$ is bounded
on $B=Y_{\theta+2\eta};$ hence $B$ is BMO-regular.  Now
$$\Phi_{B}= \frac1p\Lambda+\eta(\Phi_{X_1}-\Phi_{X_0})=
(\frac1p-\eta)\Lambda+2\eta\Phi_Z.$$ Applying Theorem 4.6
gives that $Z$ is BMO-regular.

Now we consider the converse; assume $Z$ is BMO-regular.  We
note that if $-\epsilon<\tau<1+\epsilon$ then there is a
K\"othe function space $X_{\tau}$ defined by
$\Phi_{X_{\tau}}=\Phi_{X_0}+\tau(\Phi_{X_1}-\Phi_{X_0}),$
and further each such space is $p'$-convex and $2q-$concave
where $\frac1{p'}+\frac1{2q}=1.$

We show first that if $0\le\tau_0\le 1$ then $H$ is
interpolation stable at all $0<\sigma<1$ for
$(X_{\tau_0-\frac12\epsilon},X_{\tau_0+\frac12\epsilon}).$

To this end note that $X_{\tau_0}$ is $q$-concave and so
there is a K\"othe function space $W$ defined by
$\Phi_{W}=\Phi_{X_{\tau_0}}-\frac12\epsilon\Lambda.$

Now we also have that $Z^*$ is BMO-regular.  Hence both
$Y_0=Z^{\epsilon}$ and $Y_1=(Z^*)^{\epsilon}$ are both
BMO-regular.  Now $H$ is interpolation stable at all
$0<\sigma<1$ for $(Y_0,Y_1)$ by Theorem 5.5.  Now
$\Phi_{Y_j}=\frac12\epsilon(\Lambda
+(2j-1)(\Phi_{X_1}-\Phi_{X_0}).$ Hence $$ \Phi_{Y_j}+\Phi_W
=\Phi_{X_{\tau_0}}
+(j-\frac12)\epsilon(\Phi_{X_1}-\Phi_{X_0})=\Phi_{X_{\tau_0+\epsilon
(j-\frac12)}}.$$ Thus by Lemma 5.8, $H$ is interpolation
stable at all $0<\sigma<1$ for
$(X_{\tau-\frac12\epsilon},X_{\tau_0+\frac12\epsilon}).$

Now if $I=[\alpha,\beta]$ is a closed sub-interval of
$(-\epsilon,1+\epsilon)$ we will say that $I$ is acceptable
if $H$ is interpolation stable at all $0<\sigma<1$ for
$(X_{\alpha},X_{\beta}).$ Suppose $I,J$ are two acceptable
intervals which intersect in a non-trivial interval; then we
claim that $I\cup J$ is acceptable.  In fact, excluding the
trivial cases when $I\subset J$ or $J\subset I$ we can
suppose that $I=[\alpha_1,\beta_1]$ and
$J=[\alpha_2,\beta_2]$ where
$\alpha_1<\alpha_2<\beta_1<\beta_2.$ Then we have
$H_{X_{\beta_1}}=[H_{X_{\alpha_2}},H_{X_{\beta_2}}]_{\sigma}$
where $\beta_1=(1-\sigma)\alpha_2+\sigma\beta_2$ and
similarly
$H_{X_{\alpha_2}}=[H_{X_{\alpha_1}},H_{X_{\beta_1}}]_{\sigma'}$
where $\alpha_2=(1-\sigma')\alpha_1+ \sigma'\beta_1.$ By
applying Wolff's theorem [32] we obtain
$H_{X_{\alpha_2}}=[H_{X_{\alpha_1}},H_{X_{\beta_2}}]_{\rho}$
where $\alpha_2=(1-\rho)\alpha_1+\rho\beta_2.$ It then
follows from the re-iteration theorem that we actually have
that $H$ is interpolation stable at any $0<\sigma<1$ for
$(X_{\alpha_1},X_{\beta_2})$.

Now by simple induction we can obtain that $[0,1]$ is
acceptable and this implies the result.\bull\enddemo

\demo{Remarks} Quanhua Xu has pointed out that it follows
from Theorem 5.12 that if $X_0$ is $p$-convex for some $p>1$
and if $H$ is interpolation stable at some $0<\theta<1$ then
$X_0\to X_1$ is a BMO-direction.  In fact the proof of
Theorem 5.12 essentially yields this fact, since that
direction of the argument only uses that $X_{\theta}$ is
$r$-convex, for some $r>1.$ \enddemo

\proclaim{Theorem 5.13}Suppose $H$ is a Hardy-type algebra
and that $X_0,X_1\in\Cal X.$ Suppose $X_0$ is $p$-convex for
some $p>1$ and is BMO-regular.  Suppose $X_1$ is $q$-concave
for some $q<\infty.$ Then for any $0<\theta<1$, $H$ is
interpolation stable at $\theta$ (i.e.
$[H_{X_0},H_{X_1}]_{\theta}=H_{X_{\theta}}$) where
$X_{\theta}=X_0^{1-\theta}X_1^{\theta},$ if and only if
$X_1$ is BMO-regular.\endproclaim

\demo{Proof}First note that every $X_{\theta}$ is
super-reflexive and that Proposition 3.6 can be invoked to
show the equivalence of the parenthetical statement with
interpolation stability.  One direction of the proof is
simply Theorem 5.7.  Conversely if $H$ is interpolation
stable at some $0<\theta<1$ then we may pick $0<\tau<\theta$
and $H$ is interpolation stable at $1/2$ for
$(X_{\theta-\tau},X_{\theta+\tau})$.  Hence
$(X_{\theta-\tau}\to X_{\theta+\tau})$ is a BMO-direction.
Thus if $\Phi_Z=\frac12\lambda +\tau(\Phi_{X_1}-\Phi_{X_0})$
then $Z$ is BMO-regular.  Theorem 4.6 allows us to conclude
that $X_1$ is BMO-regular.  \bull\enddemo

Let us mention at this stage that, in the case of the
standard Hardy spaces on $\bold T$ pairs $X_0,X_1$ for which
$X_0\to X_1$ is a BMO-direction can be characterized neatly
by using extended indicators.  As in [18] it is possible to
extend the indicator $\Phi_X$ to any complex $f\in L_1$ with
$|f|\in \Cal I_X\cap L\log L$ by setting
$\Phi_X(f)=\int_{\bold T}f\log x \,d\lambda$ where
$|f|=xx^*$ is the Lozanovskii factorization of $|f|$ i.e.
the unique pair $x,x^*\ge 0,$ so that supp $x,x^*=$ supp $f$
and $\|x\|_X=1,$ $\|x^*\|_{X^*}=\|f\|_1.$ The extended
$\Phi_X$ is a quasilinear map with constant $4/e$ (see Lemma
5.6 of [18]).  The following theorem follows almost directly
from Theorem 9.8 of [18].  We will not give a formal proof
here, as we plan a more detailed investigation in a
subsequent paper.

\proclaim{Theorem 5.14}Suppose $S=\bold T$ and $H=N^+$ is
the Smirnov class.  If $X_0,X_1\in\Cal X$ then $X_0\to X_1$
is a BMO-direction if and only if there is a constant $C$ so
that for any $f\in H_{\infty},$
$$|\Phi_{X_1}(f)-\Phi_{X_0}(f)| \le C\|f\|_1.$$ \endproclaim

\subheading{6.  Skew projections} We now establish some
results on ``skew'' projections.  We suppose $H$ is a closed
subalgebra of $L_{\log}$ of Hardy type (of course our
principal example of interest is the Smirnov class).  If
$w>0$ a.e. and $\log w\in L_1$ then we define $\Cal R_{w}$
to be the orthogonal projection of the weighted Hilbert
space $L_2(w)$ onto its subspace $H\cap L_2(w)=H_2(w).$

\proclaim{Theorem 6.1}Suppose $H$ is of Hardy-type.  Suppose
$X\in\Cal X$ is super-reflexive and that $0\le v,w\in L_1$
satisfy $\log v,\log w\in L_1.$ Then if $\Cal R_{v},\Cal
R_{w}$ are both bounded at $X$ then $\log v-\log w\in
BMO.$\endproclaim

\demo{Proof}Clearly by duality $\Cal R_{v}$ is also bounded
at $v^{-1}X^*$ and at $L_2(v).$ It then follows easily that
$H$ is interpolation stable at any $0<\theta<1$ for
$(L_2(v),v^{-1}X^*).$ By Theorem 5.9, $Z_1$ is BMO-regular
where for $0\le f\in L_{\infty},$ we have $$ \align
\Phi_{Z_1}(f)
&=\frac{1}{2}(\Lambda(f)-\Phi_{L_2(v)}(f)+\Phi_{v^{-1}X^*})\\
&=\frac{1}{4}(\Lambda(f)+2\Phi_{X^*}(f)-\int f\log
v\,d\lambda).  \endalign $$ By similar reasoning, $H$ is
interpolation stable at any $0<\theta<1$ for $(L_2(w),X)$
and hence $Z_2$ is BMO-regular where
$$\Phi_{Z_2}(f)=\frac{1}{4}(\Lambda(f) + 2\Phi_X(f)+\int
f\log w\,d\lambda).$$ Thus $Y=Z_1^{1/2}Z_2^{1/2}$ is
BMO-regular.  But $$ \Phi_{Y}(f)= \frac12\Lambda(f)
+\frac18\int f(\log w-\log v)d\mu.$$ Hence
$L_2((vw^{-1})^{1/4})$ is BMO-regular so that $\log v-\log
w\in BMO.$\bull\enddemo

The following theorem is suggested by a result of
Coifman-Rochberg [7] on boundedness of skew projections on
weighted $L_2-$ spaces.  We observe that although we
consider much more K\"othe spaces, our result is here
restricted to projections on Hardy subspaces; however, we
plan to investigate more general results of this type in a
forthcoming paper.

\proclaim{Theorem 6.2}Suppose $H$ is of Hardy type.  Suppose
$X_0,X_1\in\Cal X$ are super-reflexive and that $0\le v,w\in
L_1$ with $\log v,\log w\in L_1.$ Suppose that $\Cal
R_v,\Cal R_w$ are both bounded on $X_0.$ If $\Cal R_v$ is
also bounded on $X_1$ then there exists $\eta>0$ so that
$\Cal R_w$ is bounded on $X_0^{1-\theta}X_1^{\theta}$ for
$0<\theta\le \eta.$\endproclaim

\demo{Proof}Since $X_0,X_1$ are super-reflexive we may
suppose that both are $p$-convex and $q$-concave where
$\frac1p+\frac1q=1$, and $1<p\le q<\infty.$ As in Theorem
5.9 if $\epsilon=\frac1{2q}$ we can define super-reflexive
spaces $X_{\tau}$ for $-\epsilon \le \tau\le 1+\epsilon.$

Since $\Cal R_v$ is a bounded on $X_0$ and $X_1$ it is easy
to see that $H$ is interpolation stable for any $0<\theta<1$
and $(X_0,X_1).$ Thus it $X_0\to X_1$ is a BMO-direction by
Theorem 5.9 and so also is $X_{-\epsilon}\to X_1.$ Hence $H$
is interpolation stable at
$\theta=\frac{\epsilon}{1+\epsilon}$ for
$(X_{-\epsilon},X_1)$; the corresponding interpolation space
is $X_0.$

Without loss of generality we can suppose that
$d\nu=wd\lambda$ is a probability measure.  Then
$L_{\log}(\nu)=L_{\log}(\lambda)$ and so we can consider $H$
as a Dirichlet-type algebra on $(S,\nu).$ It follows from
Lemma 5.10 that since $\Cal R_w$ is bounded on $X_0$, there
exists $\eta>0$ so that $\Cal R_w$ is also bounded on
$X_{\theta}$ for all $|\theta|\le\eta.$\bull\enddemo

\vskip2truecm

\subheading{7.  The vector-valued case}

Finally let us point out an application to the vector-valued
case.  Let $S$ be a Polish space and $\mu$ be a probability
measure on $S.$ Let $X$ be a K\"othe function space on $S$
and let $Y$ be a K\"othe function space on $\bold T.$ We
denote by $Y(X)$ the K\"othe function space on $\bold
T\times S$ with measure $\lambda\times\mu$ given by
$\|f\|_{Y(X)}=\|F\|_Y$ where $F(t)=\|f(t,\cdot)\|_X.$

\proclaim{Lemma 7.1}Suppose $X\in \Cal X(S)$ and $Y\in\Cal
X(\bold T).$ Then for $0\le f\in L_{\infty}(\bold T\times
S)$ we have $$\Phi_{Y(X)}(f) = \Phi_Y(F) + \int_{\bold
T}\Phi_X(f_t)d\lambda(t)$$ where $F(t)=\int_S
f(t,s)d\mu(s),$ and $f_t(s)=f(t,s).$\endproclaim

\demo{Proof} Let us suppose first that $\|f\|_1=1$ and $f$
is a simple function of the form
$f=\sum_{j=1}^nc_j\chi_{A_j\times B_j}.$ Suppose the
Lozanovskii factorization of $F$ for $(Y,Y^*)$ is given by
$F=GH.$ Then for each $t$ suppose
$f(t,s)F(t)^{-1}=u(t,s)v(t,s)$ is the Lozanovskii
factorization for $(X,X^*)$.  Then the Lozanovskii
factorization for $(Y(X),Y(X)^*)$ is given by $f=gh$ where
$g(t,s)=G(t)u(t,s)$ and $h(t,s)=H(t)v(t,s).$ Thus $$ \align
\Phi_{Y(X)}(f) &= \int_{\bold T}\int_S f(t,s)(\log G(t)+\log
u(t,s))d\mu(s)d\lambda(t)\\ &= \Phi_Y(F) +\int_{\bold T}
\Phi_X(f_t)d\lambda(t).  \endalign $$ For general $f$ the
measurability of the integrand and the same formula follows
by a simple continuity argument (cf.  [18], Lemma
4.3).\bull\enddemo

If $X$ is a super-reflexive K\"othe function space in $\Cal
X(S)$ and $Y$ is a super-reflexive K\"othe function space on
$\bold T$ with $Y\in\Cal X(\bold T),$ then we set $H_Y(X)$
to be the closed subspace of $Y(X)$ of all functions
$f(\cdot,s)\in N^+$ for $\mu-$a.e.  $s\in S.$

We will denote the Riesz projection on $L_2(\bold T)$ by
$\Cal R$ and the vector-valued Riesz projection on
$L_2(\bold T\times S)$ by $\tilde{\Cal R}.$

We are in effect studying the Hardy-type algebra $\Cal H$
consisting of all $f\in L_{\log}(\bold T\times S)$ with
$f^s=f(\cdot,s)\in N^+$ for a.e.  $s\in S.$ For this algebra
$\Cal H_{1}$ consists of all $f\in L_{1}(\bold T\times S)$
so that $f^s\in H_{1}(\bold T)$ for a.e.  $s\in S.$ The
corresponding BMO-space we denote $\cal BMO.$

In the vector-valued case we must consider the notion of
UMD-spaces as introduced and studied initially by Burkholder
[4].  In fact a result of Bourgain [3] implies that if
$X\in\Cal X(S)$ then $X$ is a UMD-space if and only if the
Riesz projection $\tilde{\Cal R}$ is bounded on $L_2(X)$.
This characterization will be all that we require.

Now let us say that a K\"othe function space $X\in\Cal X(S)$
is {\it UMD-regular} if for some $0<\theta<1$ the space
$L_2^{\theta}X^{1-\theta}$ is a UMD-space.  If $X_0,X_1$ are
two K\"othe function spaces on $S$ we say that $X_0\to X_1$
is a {\it UMD-direction} if the space $Z$ is UMD-regular
where
$\Phi_Z=\frac12\Lambda+\frac12(\Phi_{X_1}-\Phi_{X_0}).$

\proclaim{Proposition 7.2} If $f\in L_1(\bold T\times S)$
then $f\in\cal BMO$ if and only if $f^s\in BMO$ for a.e.
$s\in S$ with $\|f\|_{\cal BMO}=\|\|f^s\|_{\cal
BMO}\|_{\infty}<\infty.$ \endproclaim

\demo{Proof}If $f\in L_1$ the map $s\to \|f^s\|_{BMO}$ is
easily seen to be measurable, and it is trivial to check
that $\|f\|_{\cal BMO}\ge \| \|f^s\|_{BMO}\|_{\infty}.$ For
the converse it is enough to note that the set $K$ of
$(\phi,\psi)$ in $L_1(\bold T)\times L_1(\bold T)$ such that
$\phi\in H_1$ and $\|\phi-\psi\|_{\infty}\le 1$ is a Borel
set.  It follows by standard selection theorems that there
is a universally measurable map $\psi\to\tilde\psi$ from
$\{\psi\in L_1;\|\psi\|_{BMO}\le 1\}$ to $H_1$ so that
$\|\psi-\tilde\psi\|_{\infty}\le 1.$ It follows easily that
$f \in L_1(\bold T\times S$ with $\|f^s\|_{BMO}\le 1$ for
a.e.  $s$ then there exists $g\in \Cal H_1$ with
$\|f-g\|_{\infty}\le 1.$\enddemo

\proclaim{Proposition 7.3}Suppose $Y\in\Cal X(\bold T)$ and
$X\in \Cal X(S).$ \newline (1) If $\tilde {\Cal R}$ is
bounded on $Y(X)$ then $\Cal R$ is bounded on $Y.$\newline
(2) If $\Cal R$ is bounded on $Y$ then $\tilde {\Cal R}$ is
bounded on $Y(L_2).$ \endproclaim

\demo{Proof}(1) Pick any fixed $0\neq x\in X$ and restrict
$\tilde{\Cal R}$ to the space $Y([x])$ where $[x]$ is the
one-dimensional space $\bold C x$.

(2) It follows directly from Krivine's theorem ([23], [24])
that the operator $(x_n)\to (\Cal Rx_n)$ is bounded on
$Y(\ell_2)$ which implies the result.\bull\enddemo

\proclaim{Proposition 7.4}Suppose $Y\in\Cal X(\bold T)$ is
super-reflexive and $X$ is a super-reflexive K\"othe
function space on $S$ with $X\in\Cal X(S).$ Then the
following conditions are equivalent:\newline (1) $Y(X)$ is a
$\tilde{\Cal R}-$direction space.  \newline (2) $Y$ is
$BMO(\bold T)-$regular and $X$ is UMD-regular.  \newline (3)
There exist constants $(C,M)$ so that if $0\le f\in Y(X)$
there exists $g\ge f$ with $\|g\|_{Y(X)}\le M\|f\|_{Y(X)}$
and $\text{ess sup}\|\log g^s\|_{BMO}\le C$ where
$g^s(t)=g(t,s)$ for $s\in S.$ \newline (4) $Y(X)$ is $\cal
BMO-$regular.

\endproclaim

\demo{Proof}Of course (3) just restates (4), and so the
equivalence of (1), (3) and (4) is just Proposition 5.11.
Let us prove that $(1)\Rightarrow (2).$ Since $Y(X)$ is a
$\tilde{\Cal R}$-direction space it follows that there
exists $\theta>0$ so that $\tilde{\Cal R}$ is bounded on
$Y_{\theta}(X_{\theta})$ where
$Y_{\theta}=L_2^{1-\theta}Y^{\theta}$ and
$X_{\theta}=L_2^{1-\theta}X^{\theta}.$ Thus by Proposition
7.3, $\Cal R$ is bounded on $Y_{\theta}$ which implies that
$Y$ is $BMO(\bold T)-$regular.  Further $\tilde{\Cal R}$ is
bounded on $Y_{\theta}(L_2)$ so that this is a $\cal
BMO-$regular space.  Hence $Y(L_2)$ is a $\cal BMO-$regular
space.  We show that $L_2(X)$ is a $\cal BMO-$regular space.
In fact if $0\le f\in L_{\infty}$ $$ \Phi_{L_2(X)}(f) =
\frac12\Lambda(F) + \int_{\bold T}\Phi_X(f_t)d\lambda$$
where $F,f_t$ are as in Lemma 7.1.  Thus $$ \Phi_{L_2(X)}(f)
= \Phi_{Y(X)}(f) -\Phi_{Y(L_2)}(f)+\Phi_{L_2(L_2)}(f)$$
whence $L_2(X)$ is $\cal BMO-$regular by Lemma 4.6.  This
implies that $\cal R$ is bounded at $L_2(X_{\phi})$ for some
$\phi>0$ i.e.  $X_{\phi}$ is UMD and so $X$ is UMD-regular.

In the converse direction we show that (2) implies that
$Y(X)$ is $\cal BMO-$regular.  Indeed if $Y$ is a
BMO-regular space then Proposition 7.3 implies that $Y(L_2)$
is $\cal BMO-$regular.  If $X$ is UMD-regular then $L_2(X)$
is a $\cal BMO$-regular space.  As in the preceding argument
we can then use Lemma 4.6 to get that $Y(X)$ is $\cal
BMO-$regular.  \bull\enddemo

\proclaim{Theorem 7.5}Suppose $(X_0,X_1)$ are
super-reflexive K\"othe function spaces in $\Cal X(S)$ and
that $(Y_0,Y_1)$ are super-reflexive K\"othe function spaces
on $\bold T$ in $\Cal X(\bold T).$ Suppose $0<\theta<1$ and
that $Y_{\theta}=[Y_0,Y_1]_{\theta}$ and
$X_{\theta}=[X_0,X_1]_{\theta}.$ Then
$[H_{Y_0}(X_0),H_{Y_1}(X_1)]_{\theta}=H_{Y_{\theta}}(X_{\theta})$
if and only if $Y_0\to Y_1$ is a BMO-direction and $X_0\to
X_1$ is a UMD-direction.\endproclaim

\demo{Proof}The necessary and sufficient condition of
Theorem 5.12 is that $Y_0(X_0)\to Y_1(X_1)$ is a $\cal
BMO-$direction.  This means by Lemma 7.1 that $W(Z)$ is a
$\cal BMO-$regular space where
$\Phi_W=\frac12(\Lambda+\Phi_{Y_1}-\Phi_{Y_0})$ and
$\Phi_Z=\frac12(\Lambda+\Phi_{X_1}-\Phi_{X_0}).$ The
equivalence of this with the fact that $W$ is BMO-regular
and $Z$ is UMD-regular is proved in Proposition 7.4.  Thus
the theorem is immediate.  \bull \enddemo

\demo{Remark}The restriction that $X_0,X_1\in\Cal X(S)$ can
easily be removed.  It is well-known that for general
K\"othe function spaces there exist weight functions $w_j,\
j=0,1$ so that $L_{\infty}\subset w_jX_j\subset L_1$.  Then
if $\tilde w_j(s,t)=w_j(s),$
$$[H_{Y_0}(X_0),H_{Y_1}(X_1)]_{\theta}=\tilde
w_0^{1-\theta}\tilde
w_1^{\theta}[H_{Y_0}(w_0X_0),H_{Y_1}(w_1X_1)]_{\theta}$$ and
this coincides with $\tilde w_0^{1-\theta}\tilde
w_1^{\theta}H_{Y_{\theta}}(w_
0^{1-\theta}w_1^{\theta}X_{\theta})$ and so on.\enddemo

We may also give a non-super-reflexive version:

\proclaim{Theorem 7.6} Suppose $X$ is a K\"othe function
spaces in $\Cal X(S)$ which is $q$-concave for some
$q<\infty$.  Suppose $(Y_0,Y_1)$ are BMO-regular K\"othe
function spaces on $\bold T$ in $\Cal X(\bold T).$ Suppose
$Y_0$ is $p$-convex where $p>1$ and that $Y_1$ is
$q$-concave.  Suppose $0<\theta<1$ and that
$Y_{\theta}=[Y_0,Y_1]_{\theta}$ and $X_{\theta}=X^{\theta}.$
Then
$[H_{Y_0}(L_{\infty}),H_{Y_1}(X)]_{\theta}=H_{Y_{\theta}}(X_{\theta})$
if and only if $X$ is UMD-regular.\endproclaim

\demo{Proof}In fact the special properties of $L_{\infty}$
imply that $Y_0(L_{\infty})$ is $\cal BMO-$regular.  Thus
from Theorem 5.12 we see that the conclusion holds if and
only if $Y_1(X)$ is $\cal BMO-$regular.  This occurs if and
only if the super-reflexive space $Y_1^{1/2}(X^{1/2})$ is
$\cal BMO-$regular or, by Proposition 7.4, if and only if
$X$ is UMD-regular.\bull\enddemo

Let us finally relate our work to that of Kisliakov and Xu
([20],[21]).  They introduce a technical condition on a
space $L_p(X,w)=w^{-1/p}L_p(X)$ where $w>0$ is a weight
function on $\bold T$ and consider when such spaces ``admit
sufficiently many analytic partitions of unity.''  Let us
say, without defining this concept precisely that $L_p(X,w)$
has the KX-property.  They show that if $X^{\alpha}$ is UMD
for some $\alpha>0$ and $\log w\in BMO$ then $L_p(X,w)$ has
the KX-property.  They also show that if $X_0,X_1$ are both
reflexive and $L_{p_0}(X_0,w_0)$ and $L_{p_1}(X_1,w_1)$ have
the KX-property then indeed $\Cal H$ is interpolation stable
for every $0<\theta<1$ for
$(L_{p_0}(X_0,w_0),L_{p_1}(X_1,w_1)).$

\proclaim{Proposition 7.7}If $X\in\Cal X(S)$ is
super-reflexive and $1<p<\infty$ is such that $L_p(X,w)$ has
the KX-property then $X$ is UMD-regular and $\log w\in BMO.$
In particular if $X^{\alpha}$ is UMD for some $\alpha>0$
then $X$ is UMD-regular.\endproclaim

\demo{Proof}As noted above, if $L_p(X,w)$ has the
KX-property then $\cal H$ is stable at any $0<\theta<1$ for
$(L_{2}(L_{2}),L_p(X,w)).$ Thus, by Theorem 7.5, $X$ is
UMD-regular and $L_p(w)$ is BMO-regular which implies that
$\log w\in BMO.$ \bull\enddemo

Note that the assumption that $X\in\Cal X(S)$ can easily be
removed by a change of weight.  Thus our results, at least
for super-reflexive spaces, include those of Kisliakov and
Xu; in fact, the conclusion also holds for spaces $X$ with
nontrivial concavity by a minor modification.

We also note that UMD-regularity of a super-reflexive
K\"othe function space is actually an isomorphic invariant;
thus if $X$ and $Y$ are two such function spaces are
isomorphic (but not necessarily as lattices), then it may be
shown that $X$ is UMD-regular if and only if $Y$ is
UMD-regular.  This can be done by the methods of [19].  Let
us conclude by remarking that in [18] we construct a
super-reflexive K\"othe function space which is not
UMD-regular.  However we do not know any example of a
super-reflexive UMD-regular space which is not already a
UMD-space (although $L_1$ is UMD-regular and not UMD).

\newpage

\subheading{References}

\item{1.} K. Barbey and H. K\"onig, {\it Abstract analytic
function theory and Hardy algebras,} \newline Springer
Lecture Notes Vol. 593, Berlin, 1977.

\item{2.} O. Blasco and Q. Xu, Interpolation between
vector-valued Hardy spaces, J. Functional Analysis 102
(1991) 331-359.

\item{3.} J. Bourgain, Some remarks on Banach spaces in
which martingale differences are unconditional, Ark.  Mat.
21 (1983) 163-168

\item{4.} D.L.  Burkholder, A geometrical characterization
of Banach spaces in which martingale difference sequences
are unconditional, Ann.  Prob. 9 (1981) 997-1011.

\item{5.} A.P.  Calder\'on, Intermediate spaces and
interpolation, the complex method, Studia Math. 24 (1964)
113-190.

\item{6.} R.R.  Coifman and R. Rochberg, Another
characterisation of BMO, Proc.  Amer.  Math.  Soc. 79 (1980)
249-254.

\item{7.} R.R.  Coifman and R. Rochberg, Projections in
weighted spaces, skew projections and inversion of Toeplitz
operators, Int.  Equns and Operator Theory, 5 (1982)
145-159.

\item{8.} M. Cotlar and C. Sadosky, On the Helson-Szeg\"o
theorem and a related class of Toeplitz kernels, Proc.
Sympos.  Pure Math. 35 (1979) 383-407.

\item{9.} M. Cotlar and C. Sadosky, On some $L^p$ versions
of the Helson-Szeg\"o theorem, 305-317 in {\it Conference on
harmonic analysis in honor of A. Zygmund,} edited by W.
Beckner, A.P.  Calder\'on, R. Fefferman and P.W.  Jones,
Wadsworth 1982.

\item{10.} M. Cotlar and C. Sadosky, The Helson-Szeg\"o
theorem in $L^p$ of the bidimensional torus, Cont.  Math.
107 (1990) 19-37.

\item{11.} M. Cwikel, J. E. McCarthy and T.H.  Wolff,
Interpolation between weighted Hardy spaces, Proc.  Amer.
Math.  Soc. 116 (1992) 381-388.

\item{12.} T.W.  Gamelin, {\it Uniform algebras,} Chelsea
Publishing Co., New York, 1984.

\item{13.} T.A.  Gillespie, Factorization in Banach function
spaces, Indag.  Math. 43 (1981) 287-300.

\item{14.} I.I.  Hirschman and R. Rochberg, Conjugate
function theory in weak$^*$-Dirichlet algebras, J.
Functional Analysis 16 (1974), 359-371.

\item{15.} H. Helson and G. Szeg\"o, A problem in prediction
theory, Ann.  Mat. Pura Appl. 51 (1960) 107-138

\item{16.} P.W.  Jones, $L^{\infty}-$estimates for the
$\bar{\partial}$ problem in the half plane, Acta Math. 150
(1983) 137-152.

\item{17.} P.W.  Jones, On interpolation between $H_1$ and
$H_{\infty}$, Springer Lecture Notes Vol. 1070 (1984)
143-151.

\item{18.} N.J.  Kalton, Differentials of complex
interpolation processes for K\"othe function spaces, Trans.
Amer.  Math.  Soc 333 (1992) 479-529.

\item{19.} N.J.  Kalton, Remarks on lattice structure in
$\ell_p$ and $L_p$ when $0<p<1,$ Israel Math.  Conf.  Series
5 (1992) 121-130.

\item{20.} S.V.  Kisliakov and Q. Xu, Sur l'interpolation
des espaces $H^p,$ C.R.  Acad.  Sci.  (Paris) 313 (1991)
249-254

\item{21.} S.V.  Kisliakov and Q. Xu, Interpolation of
weighted and vector-valued Hardy spaces, Trans.  Amer.
Math.  Soc. to appear.

\item{22.} J. Komlos, A generalization of a problem of
Steinhaus, Acta Math.  Acad.  Sci.  Hungar. 18 (1967)
217-229.

\item{23.} J.L.  Krivine, Theoremes de factorisation dans
les espaces reticules, Exposes 22-23, Seminaire
Maurey-Schwartz, Ecole Polytechnique, Paris 1973-4.

\item{24.} J. Lindenstrauss and L. Tzafriri, {\it Classical
Banach spaces, II, Function spaces,} Ergeb.  Math.
Grenzeb., vol. 97, Springer-Verlag, Berlin, Heidelberg and
New York, 1979.

\item{25.} G.Y.  Lozanovskii, On some Banach lattices,
Siberian Math.  J. 10 (1969) 419-430.

\item{26.} B. Muckenhaupt, Weighted norm inequalities for
the Hardy maximal function, Trans.  Amer.  Math.  Soc. 165
(1972) 207-226.

\item{27.} P.F.X.  M\"uller, Holomorphic martingales and
interpolation between Hardy spaces, J. d'Analyse Math. to
appear.

\item{28.} G. Pisier, Some applications of the complex
interpolation method to Banach lattices, J. d'Analyse Math.
35 (1979) 264-281.

\item{29.} G. Pisier, Interpolation between Hardy spaces and
noncommutative generalizations, Pacific J. Math. 55 (1992)
341-368.

\item{30.} J.L.  Rubio de Francia, Operators in Banach
lattices and $L^2-$inequalities, Math.  Nachr. 133 (1987)
197-209.

\item{31.} A. Torchinsky, {\it Real variable methods in
harmonic analysis,} Academic Press, 1986.

\item{32.} T.H.  Wolff, A note on interpolation spaces,
Springer Lecture Notes Vol. 908 (1982) 199-204.

\item{33.} Q. Xu, Notes on interpolation of Hardy spaces,
Ann.  Inst.  Fourier (Grenoble) 42 (1992) 875-889.

\item{34.} Q. Xu, Elementary proofs of two theorems of P.W.
Jones on interpolation between Hardy spaces, to appear.

\end